\documentclass[12pt]{article}
\usepackage{amssymb,amsmath}
\usepackage{amscd}

\usepackage{color}
\usepackage{amsthm}
\usepackage{bm}

\newcommand{\R}{\mathbb{R}}
\newcommand{\C}{\mathbb{C}}

\newtheorem{df}{Definition}[section]
\newtheorem{thm}{Theorem}[section]
\newtheorem{prop}{Proposition}[section]
\newtheorem{lem}{Lemma}[section]
\newtheorem{rem}{Remark}[section]
\newtheorem{cor}{Corollary}[section]
\newtheorem{assume}{Assumption}[section]
\begin{document}

\title{Self-similar vorticity around the boundary and non-uniqueness of solutions to the two-dimensional Navier-Stokes equations in the half space}

\author{Motofumi Aoki\footnote{Department of Mathematics, Graduate School of Science, Kyoto University, Kyoto 606-8502, Japan.  E-mail: \bf{aoki.motofumi.4e@kyoto-u.ac.jp}}
 \qquad Yasunori Maekawa\footnote{Department of Mathematics, Graduate School of Science, Kyoto University, Kyoto 606-8502, Japan. E-mail: \bf{maekawa.yasunori.3n@kyoto-u.ac.jp}}}

\date{}

\maketitle
\begin{abstract} In this paper we show the non-uniqueness of mild solutions to the two-dimensional forced Navier-Stokes equations in the half space under the noslip boundary condition, following the program established by Albritton, Bru{\'e}, and Colombo in 2022. Our construction of non-unique solutions is based on the instability of self-similar vorticity at high Reynolds numbers which concentrates around the boundary at the initial time. In our construction, therefore, a kind of boundary layer has to be taken into account in the analysis, contrasting to the known results where the unstable self-similar vorticity is located away from the boundary with $O(1)$ distance around the initial time. 
\end{abstract}

\noindent {\bf Mathematics Subject Classification (2020):35A02, 35B30, 35Q30, 35Q35, 76D05}

\section{Introduction}\label{sec.intro}

In this paper we consider the two-dimensional forced Navier-Stokes equations in the half space:
\begin{align}\label{eq.ns}
\begin{split}
& \partial_t u  -\Delta u + u \cdot \nabla u + \nabla p = F\,, \qquad t>0,~x\in \R_+^2\,,\\
& \nabla\cdot u=0,\qquad t\geq 0,~x\in \R_+^2\,,\\
& u|_{x_2=0} =0, \qquad u|_{t=0} = 0.
\end{split}
\end{align}
Here  $x=(x_1,x_2)^\top \in\R_+^2=\{x\in \R^2~|~x_2>0\}$, $u=(u_1 (t,x),u_2(t,x))^\top$ is the unknown velocity field of the fluid, $p=p(t,x)$ is the unknown pressure field, and $F=(F_1(t,x),F_2(t,x))^\top$ is a given force specified later.  We use the standard notations of derivatives; $\partial_j =\frac{\partial}{\partial x_j}$, $\nabla = (\partial_1,\partial_2)^\top$, $\Delta=\sum_{j=1}^2 \partial_j^2$, $u\cdot \nabla =\sum_{j=1}^2 u_j \partial_j$, $\nabla \cdot u = \sum_{j=1}^2 \partial_j u_j$. 

It is well known by now that if the external force $F$ is smooth and decays fast enough at infinity, then there exists a global in time smooth solution to \eqref{eq.ns} with finite energy (see, e.g., Sohr's book \cite{Sohr} and references therein). The uniqueness of such solutions in a time interval $(0,T)$ can be  also obtained, under the smallness condition on some scale critical norm, e.g., if 
\begin{align}\label{condition.u}
\sup_{0<t<T} t^\frac14 \| u(t,\cdot)\|_{L^4(\R^2_+)} <\epsilon_*
\end{align}
holds; see also Remark \ref{rem.thm.nonlinear}. Note that the condition \eqref{condition.u} is often ensured by  taking $T$ small, for only local in time uniqueness is enough in many cases. The concept of scale criticality is originated from the fact that if $(u,p)$ is a solution to \eqref{eq.ns}, then 
\begin{align}\label{scale.u}
u_\lambda (t,x) = \lambda u(\lambda^2 t, \lambda x)\,, \quad p_\lambda (t,x) = \lambda^2 p(\lambda^2 t, \lambda x)
\end{align}
is a solution to \eqref{eq.ns} with $F$ replaced by 
\begin{align}\label{scale.F}
F_\lambda (t,x) = \lambda^3 F(\lambda^2 t, \lambda x)\,.
\end{align}
Conventionally, the space-time norms which are invariant under the scalings \eqref{scale.u}, \eqref{scale.F} are called the scale critical norms for $u,p$, and $F$, respectively. As a typical example, the norm 
\begin{align}\label{scale.F.1}
\sup_{t>0} t^\frac34 \| F(t,\cdot)\|_{L^\frac43(\R^2_+)}
\end{align}
is a scale critical norm for the external force $F$, and it is not difficult to show that, if the norm of \eqref{scale.F.1} for given $F$ is small enough, then there exists a global in time unique mild solution of \eqref{eq.ns} (i.e., the solution to the integral equation) with finite energy satisfying the smallness condition such as \eqref{condition.u}. 

The smallness condition on some scale critical norms is essential in ensuring the uniqueness of solutions to the Navier-Stokes equations, which is one of the research objects of this paper. Indeed, Albritton, Bru{\'e}, and Colombo proved in their breakthrough paper \cite{ABC1} that the three-dimensional forced Navier-Stokes equations in $\R^3\times (0,T)$ admits the non-uniqueness of solutions in the Lerey-Hopf class, where the non-uniqueness is achieved around a self-similar vortex ring. The argument for the non-uniqueness in \cite{ABC1} makes use of the linear instability of the self-similar profile, as in the program introduced by Jia and {\v S}ver{\'a}k \cite{JS}. We note that the external force $F$ taken in \cite{ABC1} belongs to $BC((0,T); L^1 (\R^3)^3)$ and the non-unique solutions $u$ in \cite{ABC1} belong to $BC((0,T); L^3 (\R^3)^3)$ (here, $BC$ means {\it bounded and continuous}), and the norms of $u$ and $F$ are large in these scale critical spaces. 
In particular, it should be emphasized here that the spaces $BC((0,T); L^1 (\R^3)^3)$ and $BC((0,T); L^3 (\R^3)^3)$ rather than $BC([0,T); L^1 (\R^3)^3)$ and $BC([0,T); L^3 (\R^3)^3)$ (i.e., the continuity up to $t=0$) are essential for the non-uniqueness. Indeed, it is known by Okabe and Tsutsui \cite{OT} that for given external force in $BC([0,T); L^1 (\R^3)^3)$ there exists a local in time mild solution to the three-dimensional Navier-Stokes equations, and the uniqueness of solutions is verified in the space $BC([0,T); L^{3,\infty} (\R^3)^3)$ (here, $L^{3,\infty} (\R^3)$ is the weak $L^3$ space in $\R^3$).  
Recently, Zhan \cite{Z} proved, among others,  that the uniqueness holds even in $BC((0,T); L^{3,\infty} (\R^3)^3)$ but under the additional condition that there exists a mild solution in $BC([0,T); L^{3,\infty} (\R^3)^3)$ with the same initial data and force. The result of \cite{Z} is regarded as a kind of weak-strong uniqueness and sharply extends the result of Lions and Masmoudi \cite{LM} (see also the book of Lemari{\'e}-Rieusset \cite{LR}) on the uniqueness in $BC([0,T); L^3(\R^3)^3)$ for the three-dimensional {\it unforced} Navier-Stokes equations so that it can be applied also to the forced Navier-Stokes equations. We also note that the non-unique solutions obtained in \cite{ABC1} must be large in $BC((0,T); L^3(\R^3)^3)$, for the bilinear estimate established by Yamazaki \cite{Y} implies that the mild solutions of the three-dimensional Navier-Stokes equations must be unique if they are small in $L^\infty(0,\infty; L^{3,\infty}(\R^3)^3)$.

The result of the non-uniqueness in \cite{ABC1} is extended by the same authors in \cite{ABC2} to the three-dimensional forced Navier-Stokes equations in a smooth bounded domain under the noslip boundary condition.   
The approach of \cite{ABC2} is based on the gluing method, where the non-unique solutions in \cite{ABC1} are pasted strictly inside the domain using a cut-off, and then glued to the smooth outer solutions. In particular, the large self-similar flow $U_{ss}(t,\cdot)$, which is the source of the non-uniqueness, is located away from the boundary with $O(1)$ distance around the initial time, and the presence of the boundary plays a negligible role in the analysis. Motivated by the work of \cite{ABC1,ABC2}, in this paper we are interested in the case when the large self-similar flow $U_{ss}(t,\cdot)$ concentrates near the boundary at the initial time in the sense that the distance between the support of $U_{ss}(t,\cdot)$ and the boundary tends to zero as $t\rightarrow 0$. In this situation it is natural to expect that the effect of the boundary is more relevant, and indeed, we need to take into account a kind of the boundary layer created by the large self-similar flow. The aim of this paper is thus to develop the program of \cite{ABC1} to the case when the scale critical singularity of the force is present near the boundary, with the aid of the boundary layer analysis. To capture the essence, we focus on the case of the flat boundary for the two-dimensional problem as described in \eqref{eq.ns}. In the case of the whole space in two-dimensions, we also refer to the work by Albritton and Colombo \cite{AC}, where the non-unique Leray-Hopf solutions are presented for the forced two-dimensional Navier-Stokes equations with hypodissipation, i.e., $-\Delta$ in \eqref{eq.ns} is replaced by $(-\Delta)^\frac{\beta}{2}$ with $0<\beta<2$,  by using the approach of \cite{ABC1}.

To state the assumption and the result of this paper, let us introduce the vorticity field $\omega = {\rm rot}\, u = \partial_1 u_2-\partial_2u_1$, and then the velocity is formally recovered from the vorticity field by the Biot-Savart formula
\begin{align}\label{def.K}
K[\omega] = \nabla^\bot (-\Delta_D)^{-1} \omega\,,
\end{align}
where $\nabla^\bot =(\partial_2,-\partial_1)^\top$ and $(-\Delta_D)^{-1}g$ denotes the unique solution to the Poisson equation $-\Delta \varphi =g$ in $\R^2_+$ and $\varphi|_{x_2=0}=0$, $\lim_{|x|\rightarrow \infty} \nabla \varphi =0$. Let $\rho (x) = e^{\frac{|x|^2}{8}}$ and let $L^2_\rho (\R^2_+)$ be the weighted $L^2$ space defined as 
\begin{align}\label{def.L^2_rho}
\begin{split}
L^2_\rho (\R^2_+) & = \{g\in L^2(\R^2_+)~|~ \rho g\in L^2(\R^2_+)\}\,, \\
\langle g,\tilde g\rangle_{L^2_\rho(\R^2_+)} &= \langle \rho g , \rho \tilde g \rangle_{L^2(\R^2_+)}\,, \quad \| g \|_{L^2_\rho}  = \| \rho g\|_{L^2(\R^2_+)}\,.
\end{split}
\end{align}

\begin{assume}\label{assume}{\rm The velocity $U_E=(U_{E,1} (x), U_{E,2}(x))$ is smooth and compactly supported in $\R^2_+$ (in particular,  ${\rm dist}\, ({\rm supp}\, U_E, \partial \R^2_+)>0$) and is a solution to the stationary Euler equations in $\R^2_+$: $-U_E\cdot \nabla U_E + \nabla P_E=0$, $\nabla \cdot U_E=0$ in $\R^2_+$. Moreover, $U_E$ is linearly unstable in the sense that the operator $\Lambda_E$ in $L^2_\rho (\R^2_+)$ defined as 
\begin{align}\label{def.Lambda_E}
\begin{split}
D(\Lambda_E) & =\{\omega \in L^2_\rho (\R^2_+)~|~\nabla \cdot ( U_E \omega) \in L^2_\rho (\R^2_+)\}\\
\Lambda_E \omega & = - \nabla \cdot (U_E \omega) - K[\omega]\cdot \nabla \Omega_E
\end{split}
\end{align}
has an isolated eigenvalue $\lambda_E$ with $\Re (\lambda_E)>0$ (here, $\Re(\lambda)$ is the real part of the complex number $\lambda$). Here $\Omega_E=\partial_1 U_{E,2}-\partial_2 U_{E,1}$ is the vorticity of the Euler flow $U_E$.
}
\end{assume}

\begin{rem}\label{rem.assume}{\rm (1) The transport operator $T\omega=-\nabla\cdot (U_E\omega)$  is realized as a skew-adjoint operator in $L^2(\R^2_+)$ in virtue of $\nabla \cdot U_E=0$, and therefore, it generates a bounded and strongly continuous semigroup in $L^2(\R^2_+)$. From this fact we also have that, since the support of $U_E$ is compact, the growth bound of the semigroup generated by $T$ in the weighted space $L^2_\rho (\R^2_+)$ is also $0$.  Since the operator $S\omega=-K[\omega]\cdot \nabla \Omega_E$ is a compact operator in $L^2_\rho (\R^2_+)$, the set $\sigma (\Lambda_E) \cap \{\lambda\in \C~|~\Re (\lambda)>0\}$ consists of isolated eigenvalues with finite algebraic multiplicities, where $\sigma (A)$ is the set of the spectrum of the operator $A$; see \cite[Corollary 2.11, Proposition 2.12]{EN} for the general theory on the growth bound and the spectrum of strongly continuous semigroups.

\noindent (2) Let $O_{cpt}$ be a smooth bounded domain such that ${\rm supp}\, U_E\subset O_{cpt} \subset \overline{O_{cpt}}\subset \R^2_+$. Then the space 
\begin{align}
L_\rho^2(\R^2_+)|_{O_{cpt}} =\{\omega \in L^2_\rho (\R_+^2)~|~{\rm supp}\, \omega \subset \overline{O_{cpt}}\}
\end{align}
is invariant under the action of $\Lambda_E$, in the sense that if $\omega\in D(\Lambda_E)\cap L^2_\rho (\R_+^2)|_{O_{cpt}}$, then $\Lambda_E \omega \in L^2_\rho (\R_+^2)|_{O_{cpt}}$. It is clear that any eigenfunction of $\Lambda_E$ belongs to $D(\Lambda_E) \cap L^2_\rho (\R^2_+)|_{O_{cpt}}$. Therefore, $\lambda$ is an eigenvalue of $\Lambda_E$  if and only if so is for $\Lambda_E|_{O_{cpt}}$, where $\Lambda_E|_{O_{cpt}}$ is the restriction of $\Lambda_E$ on $L^2_\rho (\R^2_+)|_{O_{cpt}}$. Since the spectrum of $\Lambda_E|_{O_{cpt}}$ located in the right-half plane also consists of isolated eigenvalues, we have
\begin{align}
\sigma (\Lambda_E) \cap \{\lambda\in \C~|~\Re (\lambda)>0\} = \sigma (\Lambda_E|_{O_{cpt}}) \cap \{\lambda\in \C~|~\Re (\lambda)>0\}\,.
\end{align}
}
\end{rem}

\

As stated in the following theorem, we will show an example of $U_E$ satisfying Assumption \ref{assume}, whose proof will be given in Section \ref{sec.euler}.
\begin{thm}\label{thm.euler} There exists $U_E$ satisfying Assumption \ref{assume}. 
\end{thm}

\begin{rem}{\rm Let $\tilde U^{st}$ be an unstable radial flow with compact support in $\R^2$ constructed in \cite[Proposition 2.2]{ABC1}, which is obtained from Vishik's unstable radial vortices \cite{V1,V2}; see also Albritton, Bru{\'e}, Colombo, Lellis, Giri, Janisch, and Kwon \cite{ABCLGJK}.  Then the unstable flow $U_E$ in the half space, stated in Theorem \ref{thm.euler}, is obtained as $U_E=U_{E,R}(x) := \tilde U^{st} (x_1,x_2-R)$ by taking $R>1$ large enough; see Section \ref{sec.euler} for details.
}
\end{rem}

As usual, we set $L^2_\sigma (\R^2_+) =\overline{C_{0,\sigma}^\infty(\R^2_+)}^{\|\cdot \|_{L^2(\R^2_+)}}$, where $C_{0,\sigma}^\infty (\R^2_+)=\{f\in C_0^\infty (\R^2_+)^2~|~\nabla \cdot f=0 \}$. Here $C_0^\infty (\R^2_+)$ is the set of smooth and compactly supported functions in $\R^2_+$. We denote by $\mathbb{P}_\sigma$ the Helmholtz projection, i.e., the orthogonal projection from $L^2(\R^2_+)^2$ to $L^2_\sigma (\R^2_+)$. Let $\mathbb{A}=\mathbb{P}_\sigma \Delta$ be the Stokes operator in $L^2_\sigma (\R^2_+)$, i.e.,
\begin{align}
\begin{split}
D(\mathbb{A}) & = \{f\in L^2_\sigma (\R^2_+)~|~f\in W^{1,2}_0 (\R^2_+)^2 \cap W^{2,2}(\R^2_+)^2\}\,,\\
\mathbb{A} f & = \mathbb{P}_\sigma \Delta f,\qquad f\in D(\mathbb{A})\,.
\end{split}
\end{align}
Here $W^{k,p}(\R^2_+)$ and $W^{k,p}_0(\R^2_+)$, $k\in \mathbb{N}$, $1\leq p\leq \infty$, are the standard Sobolev spaces.
Then $-\mathbb{A}$ is nonnegative self-adjoint in $L^2_\sigma (\R^2_+)$, and $\mathbb{A}$ generates a strongly continuous and analytic semigroup, called the Stokes semigroup and denoted by $\{e^{t\mathbb{A}}\}_{t\geq 0}$, in $L^2_\sigma (\R^2_+)$. 

\begin{df}\label{def.mild} Let $F\in L^1_{loc} (0,\infty; L^2 (\R^2_+)^2)$. We say that $u$ is a mild solution of \eqref{eq.ns} if $u$ satisfies the following conditions:

\noindent {\rm (i)} $u \in L^\infty (0,T; L^2_\sigma (\R_+^2))\cap C((0,T); L^2_\sigma (\R^2_+))$ and $u\in L^2 (\delta, T; W^{1,2}_0(\R^2_+)^2)$ for any $0<\delta<T<\infty$.

\noindent {\rm (ii)} $\lim_{t\rightarrow 0} u(t) =0$ weakly in $L^2_\sigma (\R^2_+)$.

\noindent {\rm (iii)} For any $0<s<t<\infty$, the equality
\begin{align}
u(t) = e^{(t-s)\mathbb{A}} u(s) - \int_s^t e^{(t-\tau)\mathbb{A}} \mathbb{P}_\sigma \nabla \cdot (u\otimes u) (\tau) d\tau + \int_s^t e^{(t-\tau)\Delta} \mathbb{P}_\sigma F(\tau) d\tau
\end{align}
holds.
\end{df}

The main result of this paper is the non-uniqueness of mild solutions stated as follows.
\begin{thm}\label{thm.nonlinear} Let $F=F_\alpha = \alpha (\partial_t u_E-\Delta u_E)$, where $\alpha$ is a positive constant and $u_E = \frac{1}{\sqrt{t}} U_E (\frac{x}{\sqrt{t}})$, and $U_E$ satisfies Assumption \ref{assume}. Then, there exists $\alpha_E\geq 1$ such that if $\alpha\geq \alpha_E$ then there exist at least two mild solutions of \eqref{eq.ns}.
\end{thm}

\begin{rem}\label{rem.thm.nonlinear}{\rm (1) Clearly, the self-similar flow $\alpha u_E$ itself is a mild solution of \eqref{eq.ns}. We note that the support of $\alpha u_E(t,\cdot)$ in $\R^2_+$ converges to the origin as $t\rightarrow 0$. 

\noindent (2) In Theorem \ref{thm.nonlinear}, taking $\alpha$ large is essential in proving the non-uniqueness. Indeed, if $\alpha$ is small enough, then we have the smallness such as 
$$\displaystyle \sup_{t>0} t^\frac34 \| F(t,\cdot)\|_{L^\frac43(\R^2_+)}\leq C\alpha\ll 1\,.$$ 
Then the semigroup estimates such as $\| e^{t\mathbb{A}} \mathbb{P}_\sigma \nabla \cdot V\|_{L^4(\R^2_+)}\leq Ct^{-\frac34} \| V\|_{L^2(\R^2_+)}$ for $V\in L^2(\R^2_+)^{2\times 2}$ and $\| e^{t\mathbb{A}} \mathbb{P}_\sigma F\|_{L^4(\R^2_+)}\leq Ct^{-\frac12} \| F\|_{L^\frac43(\R^2_+)}$ for $F\in L^\frac43(\R^2_+)^{2}$ enable us to construct the unique mild solution satisfying the smallness condition \eqref{condition.u}. That is, $\alpha u_E$ is the only solution of \eqref{eq.ns} satisfying \eqref{condition.u} if $\alpha$ is small enough.
}
\end{rem}

To prove Theorem \ref{thm.nonlinear} we basically follow the approach of \cite{ABC1}, and hence, it is useful to introduce the self-similar variables
\begin{align}
v(\tau,\xi)  = e^\frac{\tau}{2} u (e^\tau, e^\frac{\tau}{2} \xi)\,, \quad \tau\in \R\,, ~\xi=(\xi_1,\xi_2)\in \R^2_+\,.
\end{align}
Then $v$ satisfies the equations 
\begin{align*}
\partial_\tau v - (\Delta +\frac{\xi}{2}\cdot \nabla + \frac12) v + v\cdot \nabla v + \nabla q = G\,, \quad \tau\in \R\,,~\xi\in \R^2_+
\end{align*}
with the divergence free condition $\nabla \cdot v=0$ and the no-slip boundary condition $v|_{\xi_2=0}=0$. Here $G(\tau,\xi) = e^{\frac32\tau} F(e^\tau,e^\frac{\tau}{2}\xi)$, and $\Delta$ and $\nabla$ are now the differential operators about the $\xi$ variables. When $F$ is chosen as in Theorem \ref{thm.nonlinear}, the stationary flow $v(\tau,\xi)=\alpha U_E(\xi)$ clearly satisfies the equations. The other solution is then constructed around $\alpha U_E$, and therefore,  we are interested in the linearized operator in the self-similar variables $\xi$ around the profile $\alpha U_E$. This motivates us to study the operator defined as follows.
\begin{align}\label{def.L_alpha}
\begin{split}
D(\mathbb{L}_\alpha) & = \{ f\in L^2_\sigma (\R^2_+)~|~f\in W^{1,2}_0(\R^2_+)^2 \cap W^{2,2}(\R^2_+)^2,\quad \xi \cdot \nabla f\in L^2(\R^2_+)^2 \}\,,\\
\mathbb{L}_\alpha f & = \frac{1}{\alpha} \big ( \mathbb{A} + \frac{\xi}{2}\cdot \nabla +\frac12 \big ) f -  \mathbb{P}_\sigma \big (U_E \cdot \nabla f + f\cdot \nabla U_E \big )\,, \quad f\in D(\mathbb{L}_\alpha)\,.
\end{split}
\end{align}
Note that the parameter $\alpha$ originally describes the size of the external force, while it plays a role of the viscosity coefficient in \eqref{def.L_alpha}. 

Since $\lambda_E$ is an isolated unstable eigenvalue of $\Lambda_E$, there exists $r_E\in (0,\Re (\lambda_E))$ such that $\overline{B_{r_E} (\lambda_E)} \setminus \{\lambda_E\} \subset \rho_{re} (\Lambda_E)$. Here $B_{r} (\zeta) =\{\lambda\in \C~|~|\lambda-\zeta|<r\}$ and $\rho_{re} (A)$ is the resolvent set of the operator $A$. The following theorem about the existence of unstable eigenvalues of the operator $\mathbb{L}_\alpha$ is the key to show the non-uniqueness of the nonlinear problem.
\begin{thm}\label{thm.spectrum} For any $\epsilon\in (0, \frac{r_E}{2})$ there exists $\alpha_\epsilon\geq 1$ such that if $\alpha\geq \alpha_\epsilon$ then there exists an isolated eigenvalue $\lambda_\alpha$ of $\mathbb{L}_\alpha$ satisfying $|\lambda_\alpha - \lambda_E|<\epsilon$. In particular, $\mathbb{L}_\alpha$ has an unstable eigenvalue $\lambda_\alpha$. 
\end{thm}

Once Theorem \ref{thm.spectrum} is proved, Theorem \ref{thm.nonlinear} follows from the similar argument as in \cite[Section 4]{ABC1}. As for the proof of Theorem \ref{thm.spectrum}, by considering the eigenprojection around the eigenvalue $\lambda_E$ of $\Lambda_E$, the analysis is essentially reduced to the construction of the solution to the resolvent problem 
\begin{align}\label{eq.resolvent}
\lambda v - \mathbb{L}_\alpha v =f
\end{align}
for the complex number $\lambda$ belonging to $\rho_{re}(\Lambda_E)$ and for given $f$ such that $g={\rm rot}\, f\in L^2(\R^2_+)$ and ${\rm supp}\, g$ is compact in $\R^2_+$. From the condition $\lambda\in \rho_{re}(\Lambda_E)$ and $\alpha\gg 1$ it is natural to construct the solution around $K[(\lambda I - \Lambda_E)^{-1} g]$, which is the solution to the linearized Euler equations. This can be done as in the whole space case if the boundary condition is slip-type, however, it is not so straightforward if the boundary condition is noslip, due to the discrepancy between the boundary conditions in the linearized Euler equations and in the linearized Navier-Stokes equations. To handle this issue we follow the approach used in the stability analysis for the Prandtl boundary layer (e.g., see G{\'e}rard-Varet, Maekawa, and Masmoudi \cite{GVMM} and references therein), where we first solve the problem under the artificial boundary condition such as the perfect slip boundary condition (i.e., the vorticity is zero on $\partial\R^2_+$) by analyzing the vorticity equations, and the noslip boundary condition is then recovered by constructing a suitable boundary layer. In constructing the boundary layer, we make use of the advantage that the support of the profile $U_E(\xi)$ is strictly away from the boundary in the self-similar variables. In particular, the boundary layer can be taken independently of the equation itself; see the proof of Proposition \ref{prop.linear}. The instability in the high frequency for the boundary layer does not appear here in virtue of the fact that the support of the profile $U_E(\xi)$ is strictly away from the boundary, where the observation of this kind is motivated from the work of Maekawa \cite{M} that verifies the Prandtl boundary layer expansion in the inviscid limit when the initial vorticity is located away from the boundary.

This paper is organized as follows. In Section \ref{sec.pre} we collect some basic estimates related with the Biot-Savart formula \eqref{def.K} and the Stokes operator $\mathbb{A}$, where the proofs are more or less standard and are postponed to the appendix. In Section \ref{sec.slip} we construct the solution to the resolvent equations for the linearized problem but under the perfect slip boundary condition. Then the solution to \eqref{eq.resolvent} (with noslip boundary condition) is constructed in Section \ref{sec.noslip}, which leads to Theorem \ref{thm.spectrum}. Theorem \ref{thm.nonlinear} is proved in Section \ref{sec.nonlinear} as an application of the results in Section \ref{sec.noslip}.
Finally, Section \ref{sec.euler} is devoted to the proof of Theorem \ref{thm.euler}.

\section{Preliminaries}\label{sec.pre}

\subsection{Estimates for the Biot-Savart formula}\label{subsec.BS}

In this subsection we collect the estimates of the velocity $K[\omega]=(K_1[\omega],K_2[\omega])=\nabla^\bot (-\Delta_D)^{-1} \omega$ for $\omega\in L^2_\rho (\R^2_+)$.

\begin{lem}\label{lem.bs} Let $\omega\in L^2_\rho (\R^2_+)$. Then the following estimates hold.
\begin{align}\label{est.lem.bs.1}
\begin{split}
\| K[\omega]\|_{L^2(\R^2_+)}  & \leq C\| \xi_2 \omega\|_{L^2(\R^2_+)}\,,\\
\| \nabla K[\omega] \|_{L^2(\R^2_+)} & \leq \| \omega\|_{L^2(\R^2_+)}\,,\\
\| K_1[\omega]\big|_{\xi_2=0} \|_{L^2(\R)} & \leq C\| \langle \xi_2\rangle \omega\|_{L^2 (\R^2_+)}\,.
\end{split}
\end{align}
Moreover, for any $\kappa>0$ there exists $C_\kappa>0$ such that 
\begin{align}\label{est.lem.bs.2}
\| \partial_1^j K_1[\omega]\big|_{\xi_2=0} \|_{L^2(\R)} & \leq C_\kappa \big ( \| \langle \xi_2\rangle \omega\|_{L^2(\R^2_+)} + \| \partial_1^j \omega\|_{L^2 (\R\times (0,\kappa))}\big ) \,, \quad 1\leq j\leq 5\,,
\end{align}
and 
\begin{align}
\begin{split}
& \| \xi_1\partial_1^{j+1} K_1[\omega]\big|_{\xi_2=0}\|_{L^2(\R)} \\
& \leq C_\kappa \big ( \|\langle \xi\rangle \omega\|_{L^2(\R^2_+)} + \| \partial_1^j \omega\|_{L^2(\R\times (0,\kappa))} + \sum_{l=0}^j \|  \xi_1\partial_1^{1+l} \omega \|_{L^2(\R\times (0,\kappa))}\big )\,, \quad j=0,1\,.  \label{est.lem.bs.3}
\end{split}
\end{align}
Here $\langle \xi_2\rangle = (1+\xi_2^2)^\frac12$ and $\langle \xi\rangle=(1+|\xi|^2)^\frac12$.
\end{lem}

The proof of Lemma \ref{lem.bs} is postponed to the appendix.

\subsection{Estimates for the harmonic oscillator}\label{subsec.harmonic}

Let us recall that $\rho (\xi) = e^{\frac{|\xi|^2}{8}}$. 
In this subsection we recall the property of the operator $H$ in $\R^2$ defined as 
\begin{align}
\begin{split}
D(H) & = \{ g\in L^2_\rho (\R^2)~|~\Delta g\in L^2_\rho (\R^2),~\xi \cdot \nabla g\in L^2_\rho (\R^2)\}\,,\\
H g & = (\Delta  + \frac{\xi}{2} \cdot \nabla +1) g\,, \quad g\in D(H)\,.
\end{split}
\end{align}
It is well known that $-H$ is a nonnegative self-adjoint operator in $L^2_\rho (\R^2)$ with compact resolvent, and we have the identity $\rho (-H) \frac{1}{\rho} =-\Delta +\frac{|\xi|^2}{16}-\frac12$, where the right-hand side is the Hamiltonian of the harmonic oscillator. Moreover, the space $L^2_{\rho,0}(\R^2) = \{ g \in L^2_\rho (\R^2)~|~\int_{\R^2} g (\xi) d\xi =0\}$ is an invariant space under the action of $H$ and $-H\geq \frac12$ in $L^2_{\rho,0} (\R^2)$; see, e.g., \cite[Lemma 4.7]{GW}. Hence we have for $w\in L^2_{\rho,0}(\R^2)\cap D(H)$,
\begin{align}\label{est.harmonic}
\Re \langle -H w,w\rangle_{L^2_\rho (\R^2)} & = \frac18 \Re \langle -H w,w\rangle_{L^2_\rho (\R^2)} +  \frac18 \Re \langle -H w,w\rangle_{L^2_\rho (\R^2)} + \frac34 \Re \langle -H w,w\rangle_{L^2_\rho (\R^2)} \nonumber \\
& \geq \frac18 (\| \nabla w\|_{L^2_\rho (\R^2)}^2 -\| w\|_{L^2_\rho (\R^2)}^2) \nonumber \\
& \quad +  \frac18 \Re \langle (-\Delta +\frac{|\xi|^2}{16}-\frac12) (\rho w), \rho w\rangle_{L^2 (\R^2)} + \frac38 \| w\|_{L^2_\rho (\R^2)}^2\nonumber \\
& \geq \frac18 \| \nabla w\|_{L^2_\rho (\R^2)}^2 + \frac{1}{128} \| \xi w\|_{L^2_\rho(\R^2)}^2 +\frac{3}{16} \| w\|_{L^2_\rho (\R^2)}^2\,.
\end{align}
This lower bound will be used later, though the explicit constants are not essential in this paper.

\subsection{Estimates for the Stokes system in the self-similar variables}\label{subsec.stokes}

Let $\epsilon\in (0, \frac{r_E}{2})$. The parameter $\alpha$ is always taken large enough so that $\alpha \geq \frac{1}{r_E}$. We see that $\partial B_\epsilon (\lambda_E)=\{\lambda\in \C~|~|\lambda-\lambda_E|=\epsilon\}$ is contained in the resolvent set of $\mathbb{H}_\alpha = \frac{1}{\alpha} (\mathbb{A} +\frac{\xi}{2}\cdot \nabla +\frac12)$ in $L^2_\sigma (\R^2_+)$, where the domain of $\mathbb{H}_\alpha$ is given as
\begin{align*}
D(\mathbb{H}_\alpha) & = \{f\in L^2_\sigma (\R^2_+) \cap W^{2,2}(\R^2_+)^2\cap W^{1,2}_0(\R^2_+)^2 ~|~\xi \cdot \nabla f\in L^2_\sigma(\R^2_+)\}\\
& =\{f\in L^2_\sigma (\R^2_+) \cap W^{2,2}(\R^2_+)^2\cap W^{1,2}_0(\R^2_+)^2 ~|~\xi \cdot \nabla f\in L^2(\R^2_+)^2\}\,.
\end{align*}
Note that the second equality holds due to the invariance: if $f\in L^2_\sigma (\R^2_+) \cap W^{2,2}(\R^2_+)^2\cap W^{1,2}_0(\R^2_+)^2$ and $\xi \cdot \nabla f\in L^2(\R^2_+)^2$, then $\xi \cdot \nabla f\in L^2_\sigma (\R^2_+)$. The estimates of $v = \big (\lambda I -\mathbb{H}_\alpha \big )^{-1} f$ for $f\in L^2_\sigma (\R^2_+)$ and $\lambda\in \partial B_\epsilon (\lambda_E)$ are stated as follows. 

\begin{lem}\label{lem.stokes} Under the definition above, we have for $v=(\lambda I - \mathbb{H}_\alpha)^{-1} f$, 
\begin{align}\label{est.lem.stokes.1}
\| v\|_{L^2(\R^2_+)} + \frac{1}{\alpha^\frac12} \| \nabla v \|_{L^2(\R^2_+)} + \frac{1}{\alpha} \big ( \| \nabla^2 v \|_{L^2(\R^2_+)} + \| \xi\cdot \nabla v\|_{L^2(\R^2_+)}\big ) \leq C \| f\|_{L^2(\R^2_+)}\,.
\end{align}
Here the constant $C$ is taken uniformly in $\lambda\in \partial B_\epsilon (\lambda_E)$ and $\alpha \geq \frac{1}{r_E}$. Moreover, for smooth bounded domains $O_{cpt}, O_{cpt}'$ satisfying $\overline{O_{cpt}}\subset O_{cpt}'\subset \overline{O_{cpt}'} \subset \R^2_+$, we have 
\begin{align}\label{est.lem.stokes.2}
\| \nabla^k {\rm rot} \, v \|_{L^2(O_{cpt})} \leq C \big ( \| f\|_{L^2(\R^2_+)} + \sum_{j=0}^k \| \nabla^j {\rm rot}\, f\|_{L^2(O_{cpt}')}\big )\,,\quad k=0,1\,.
\end{align}
Here $C$ depends only on $O_{cpt}$ and $O_{cpt}'$, and in particular, $C$ is independent of $\lambda\in \partial B_\epsilon (\lambda_E)$ and $\alpha \geq \frac{1}{r_E}$.
\end{lem}

The proof of Lemma \ref{lem.stokes} is postponed to the appendix.
Let $\mathbb{H}=\alpha \mathbb{H}_\alpha= \mathbb{A}+\frac{\xi}{2}\cdot \nabla + \frac12$, and let $\{e^{\tau\mathbb{H}}\}_{\tau\geq 0}$ be the semigroup generated by $\mathbb{H}$ in $L^2_\sigma (\R^2_+)$, which is strongly continuous in $L^2_\sigma (\R^2_+)$. We observe that $v(\tau,\xi) = (e^{\tau \mathbb{H}}f)(\xi)$ is related to the original Stokes semigroup $\{e^{t\mathbb{A}}\}_{t\geq 0}$ by the scaling 
\begin{align}\label{proof.lem.stokes.semigroup.1}
v(\tau,\xi) = e^\frac{\tau}{2} u(e^\tau-1,e^\frac{\tau}{2}\xi)\,,\quad  u(t,x) = (e^{t\mathbb{A}}f)(x)\,. 
\end{align}
For the Stokes semigroup $\{e^{t\mathbb{A}}\}_{t\geq 0}$ (see, e.g., \cite{S,U,Sohr}) we have 
\begin{align}\label{proof.lem.stokes.semigroup.2}
\| e^{t\mathbb{A}} f\|_{L^2(\R^2_+)} + t^\frac12 \| \nabla e^{t\mathbb{A}} f\|_{L^2(\R^2_+)} + t^\frac12 \| e^{t\mathbb{A}} f\|_{L^\infty(\R^2_+)} \leq C \| f\|_{L^2(\R^2_+)}\,.
\end{align}
Indeed, we note that the estimate $\| e^{t\mathbb{A}} f\|_{L^2(\R^2_+)} + t^\frac12 \| \nabla e^{t\mathbb{A}} f\|_{L^2(\R^2_+)} \leq C\| f\|_{L^2(\R^2_+)}$ follows from the standard $L^2$ theory of the Stokes operator, while the estimate $t^\frac12 \| e^{t\mathbb{A}} f\|_{L^\infty(\R^2_+)} \leq C \| f\|_{L^2(\R^2_+)}$ follows from the embedding $\| e^{t\mathbb{A}} f\|_{L^\infty(\R^2_+)} \leq C \| \nabla e^{t\mathbb{A}} f\|_{L^4(\R^2_+)}^\frac12 \| e^{t\mathbb{A}} f\|_{L^4(\R^2_+)}^\frac12$ and the $L^p$-$L^q$ estimates stated in Ukai \cite[Theorem 3.1]{U}. We also have from Koboyashi and Kubo \cite{KK} the following weighted estimate
\begin{align}\label{proof.lem.stokes.semigroup.3}
\| \langle x\rangle^m e^{t\mathbb{A}} f\|_{L^2(\R^2_+)} \leq C (1+t^\frac{m}{2}) \| \langle x\rangle^m f\|_{L^2(\R^2_+)}\,,\qquad t>0,
\end{align}
for $0<m<\frac12$ and $f\in L^2_\sigma (\R^2_+)$ with $\langle x\rangle^m f\in L^2(\R^2_+)^2$, where $\langle x\rangle = (1+|x|^2)^\frac12$. Thus, the identity \eqref{proof.lem.stokes.semigroup.1} and the estimates \eqref{proof.lem.stokes.semigroup.2}-\eqref{proof.lem.stokes.semigroup.3} give the following estimates for $\{e^{\tau \mathbb{H}}\}_{\tau\geq 0}$.
\begin{lem}\label{lem.stokes.semigroup} There exists $C>0$ such that, for any $f\in L^2_\sigma (\R^2_+)$ and $\tau>0$, 
\begin{align}\label{est.lem.stokes.semigroup}
\| e^{\tau \mathbb{H}} f\|_{L^2(\R^2_+)} + a(\tau)^\frac12 \| \nabla e^{\tau \mathbb{H}} f\|_{L^2(\R^2_+)} + a(\tau)^\frac12 \| e^{\tau \mathbb{H}} f\|_{L^\infty (\R_+^2)} \leq C \|f\|_{L^2(\R^2_+)}\,.
\end{align}
Here $a(\tau) = 1-e^{-\tau}$. Moreover, if $0<m<\frac12$ and $\langle \xi\rangle^m f\in L^2(\R^2_+)^2$ in addition, then we have
\begin{align}\label{est.lem.stokes.semigroup'}
\| \langle e^\frac{\tau}{2} \xi\rangle^m e^{\tau \mathbb{H}} f\|_{L^2(\R^2_+)} \leq C e^\frac{m\tau}{2} \| \langle \xi\rangle^m f\|_{L^2(\R^2_+)}\,.
\end{align}

\end{lem}

\section{Linear analysis under perfect slip boundary condition}\label{sec.slip}

In this section we consider the resolvent problem 
\begin{align}\label{eq.slip}
\begin{split}
& \lambda v  - \frac{1}{\alpha} (\Delta +\frac{\xi}{2}\cdot \nabla +\frac12) v + U_E \cdot \nabla v + v\cdot \nabla U_E + \nabla p = f\,, \quad \xi \in \R_+^2\,,\\
& \nabla\cdot v=0,\quad \xi \in \R_+^2\,,\\
& \omega|_{\xi_2=0} = v_2|_{\xi_2=0} =0\,.
\end{split}
\end{align}
Here $\lambda\in \C$ with $\Re (\lambda)>0$ and $\omega = {\rm rot}\, v=\partial_1v_2-\partial_2 v_1$ is the vorticity field of $v$. Assuming that $g:={\rm rot}\, f\in L^2_\rho (\R^2_+)$ in addition, by taking ${\rm  rot}$ in the first equation, the system \eqref{eq.slip} is equivalent with the linearized vorticity equation 
\begin{align}\label{eq.slip.vor}
\begin{split}
(\lambda I  - M_\alpha ) \omega  = g\,,
\end{split}
\end{align}
where the operator $M_\alpha$ in $L^2_\rho(\R^2_+)$ is defined as 
\begin{align}
\begin{split}
D(M_\alpha) & = \{ \omega\in L^2_\rho (\R^2_+)~|~\Delta \omega \in L^2_\rho (\R^2_+)\,,~\xi \cdot \nabla \omega\in L^2_\rho (\R^2_+)\,,~ \omega|_{\xi_2=0}=0\}\,, \\
M_\alpha \omega & = \frac{1}{\alpha} (\Delta +\frac{\xi}{2}\cdot \nabla +1) \omega +\Lambda_E \omega\,, ~\omega \in D(M_\alpha)\,,
\end{split}
\end{align}
and $\Lambda_E \omega = - U_E\cdot \nabla \omega - K[\omega]\cdot \nabla \Omega_E$ with $K[\omega]=\nabla^\bot (-\Delta_D)^{-1}\omega$ for $\omega\in D(M_\alpha)$. By extending $\omega$ and $U_E$ to $\R^2$ using the reflection 
\begin{align*}
\tilde \omega (\xi) & = -\omega (\xi_1,-\xi_2),\qquad \xi_2<0\,,\\
\tilde U_{E,1} (\xi) & = U_{E,1} (\xi_1,-\xi_2),\qquad \xi_2<0\,,\\
\tilde U_{E,2} (\xi) & = - U_{E,2} (\xi_1,-\xi_2),\qquad \xi_2<0\,,
\end{align*}
the system \eqref{eq.slip.vor} is reduced to the system in $\R^2$:
\begin{align}\label{eq.slip.vor.tilde}
\begin{split}
& \lambda \tilde \omega  - \frac{1}{\alpha} H \tilde \omega + \tilde U_E \cdot \nabla \tilde \omega + \tilde K[\tilde \omega]\cdot \nabla \tilde \Omega_E  = \tilde g\,, \quad \xi \in \R^2\,.
\end{split}
\end{align}
Here $H=\Delta +\frac{\xi}{2}\cdot \nabla +1$ is the operator discussed in Subsection \ref{subsec.harmonic},  $\tilde K[\tilde \omega] = \nabla^\bot (-\Delta)^{-1}\tilde \omega$ is the Biot-Savart formula in $\R^2$, and $\tilde g$ is the odd extension of $g$. In this section we consider the system \eqref{eq.slip.vor.tilde} in the weighted space $L^2_{\rho,odd} (\R^2)$, where  
 \begin{align}\label{def.L^2_rho.odd}
 \begin{split}
 L^2_{\rho,odd} (\R^2) & = \{\tilde g \in L^2 (\R^2)~|~\rho \tilde g\in L^2(\R^2),~\tilde g(\xi_1,-\xi_2) = -\tilde g(\xi_1,\xi_2)\}\,,\\
\| \tilde g\|_{L^2_{\rho,odd} (\R^2)} & =  \| \rho \tilde g\|_{L^2 (\R^2)}\,.
\end{split}
\end{align}
Then we denote by $\tilde H$ the restriction of $H$ on $L^2_{\rho,ood}(\R^2)$, i.e., 
\begin{align}\label{def.tilde.H}
\begin{split}
D(\tilde H) & = \{\tilde g\in L^2_{\rho,odd}(\R^2)~|~\Delta \tilde g\in L^2_{\rho,odd}(\R^2), ~\xi\cdot \nabla \tilde g\in L^2_{\rho,odd} (\R^2)\}\,,\\
\tilde H \tilde g & = H \tilde g\,, \quad g\in D(\tilde H)\,,
\end{split}
\end{align}
and we also set
\begin{align}\label{def.tilde.Lambda_E}
\begin{split}
D(\tilde \Lambda_E) & = \{\tilde g\in L^2_{\rho,odd}(\R^2)~|~\nabla \cdot (\tilde U_E \tilde g)\in L^2_{\rho,odd}(\R^2)\}\,,\\
\tilde \Lambda_E \tilde g & = -\nabla \cdot (\tilde U_E \tilde g) - \tilde K[\tilde g] \cdot \nabla \tilde \Omega_E\,, \quad \tilde g\in D(\tilde \Lambda_E)\,.
\end{split}
\end{align}
By Assumption \ref{assume}, the extended operator $\tilde \Lambda_E$ also has an isolated eigenvalue $\lambda_E$ with $\Re (\lambda_E)>0$. Let us set
\begin{align}
D(\tilde M_\alpha) = D(\tilde H),\quad \tilde M_\alpha = \frac{1}{\alpha} \tilde H + \tilde \Lambda_E.
\end{align}
Then \eqref{eq.slip.vor.tilde} is written as $(\lambda I - \tilde M_\alpha ) \tilde \omega = \tilde g$ in $L^2_{\rho,odd}(\R^2)$. Let us recall that $\rho_{re} (\tilde M_\alpha)$ denotes the resolvent set of $\tilde M_\alpha$. 

\begin{thm}\label{thm.slip.tilde} For any $\epsilon \in (0,\frac{r_E}{2})$ there exists $\alpha_\epsilon'\geq 1$ such that if $\alpha\geq \alpha_\epsilon'$ then $\partial B_\epsilon (\lambda_E) \subset \rho_{re} (\tilde M_\alpha)$ and 
\begin{align}\label{est.thm.slip.1}
& \sup_{\lambda\in \partial B_\epsilon (\lambda_E),~\alpha\geq \alpha_\epsilon'} \| (\lambda I -\tilde M_\alpha)^{-1}\|_{L^2_{\rho,odd} (\R^2)\rightarrow L^2_{\rho,odd} (\R^2)} <\infty\,,
\end{align}
and for any  $\tilde g\in L^2_{\rho,odd} (\R^2)$, 
\begin{align}\label{est.thm.slip.2}
\lim_{\alpha\rightarrow \infty} \sup_{\lambda\in \partial B_\epsilon (\lambda_E)} \|(\lambda I - \tilde M_\alpha)^{-1} \tilde g - (\lambda I - \tilde \Lambda_E)^{-1} \tilde g \|_{L^2_{\rho,odd} (\R^2)}=0\,. 
\end{align}
\end{thm}

\begin{proof} Since we are working in the whole space, the argument of the proof is similar to the one of \cite[Lemma 3.4]{ABC1}.  To simplify the notation we write $\|\cdot\|_{L^2_\rho (\R^2)}$ for $\| \cdot \|_{L^2_{\rho,odd}(\R^2)}$. Let us decompose $\tilde M_\alpha$ as 
\begin{align*}
\tilde M_\alpha  = \tilde T_\alpha + \tilde S\,,\quad \tilde T_\alpha \tilde g & = \frac{1}{\alpha} \tilde H \tilde g - \tilde U_E \cdot \nabla \tilde g \,,\quad \tilde S\tilde g = -\tilde K[\tilde g] \cdot \nabla \Omega_E\,.
\end{align*}
We note that 
\begin{align*}
\Re \langle (\lambda - \tilde T_\alpha) w,w\rangle_{L^2(\R^2)} = (\Re (\lambda) -\frac{1}{2\alpha}) \|w\|_{L^2(\R^2)}^2 + \frac{1}{\alpha}\| \nabla w\|_{L^2(\R^2)}^2
\end{align*}
for any $w\in D(\tilde T_\alpha)=D(\tilde H)$. 
If $\alpha \geq \frac{2}{r_E}$, $r_E=\Re (\lambda_E)$, and $\lambda\in \partial B_\epsilon (\lambda_E)$ with $\epsilon\in (0,\frac{r_E}{2})$, then $\Re (\lambda)-\frac{1}{2\alpha}\geq \frac{r_E}{2}-\frac{1}{2\alpha} \geq \frac{r_E}{4}$, and therefore, $\lambda I - \tilde T_\alpha$ is injective in $L^2_{odd} (\R^2) \subset L^2(\R^2)$ (with the trivial definition of $L^2_{odd}(\R^2)$). Hence, $\lambda I - \tilde T_\alpha$ is injective in the weighted space $L^2_{\rho,odd} (\R^2)$ as well, which implies that $\lambda I -\tilde T_\alpha$ is invertible in $L^2_{\rho,odd} (\R^2)$, since the spectrum of $\tilde T_\alpha$ realized in $L^2_{\rho,odd}(\R^2)$ consists only of (isolated) eigenvalues. Then, the lower bound for the above bilinear form implies 
\begin{align}\label{proof.thm.slip.1}
\| (\lambda I - \tilde T_\alpha)^{-1} \tilde g\|_{L^2(\R^2)} \leq \frac{4}{r_E} \| \tilde g\|_{L^2(\R^2)},\quad \tilde g\in L^2_{\rho,odd} (\R^2)\,.
\end{align} 
Next we see for $w\in D(\tilde H)$,
\begin{align*}
& \Re \langle (\lambda - \tilde T_\alpha) w,w\rangle_{L^2_\rho (\R^2)} \\
& = \Re (\lambda) \| w\|_{L^2_\rho (\R^2)}^2 + \frac{1}{\alpha} \Re \langle  -H w,w\rangle_{L^2_\rho (\R^2)} + \Re \langle \tilde U_E\cdot \nabla w, w\rangle_{L^2_\rho (\R^2)}\,,
\end{align*}
which implies from \eqref{est.harmonic} and from the integration by parts for the term $\Re \langle \tilde U_E\cdot \nabla w, w\rangle_{L^2_\rho (\R^2)}$ that 
\begin{align*}
\Re (\lambda) \| w\|_{L^2_\rho (\R^2)}^2  + \frac{1}{\alpha} \| \nabla w \|_{L^2_\rho (\R^2)}^2  \leq C( \frac{1}{\Re (\lambda)} \| (\lambda - \tilde T_\alpha ) w\|_{L^2_\rho (\R^2)}^2 + \| w\|_{L^2(\R^2)}^2 )\,.
\end{align*}
Here the fact that the support of $U_E$ is compact is used to derive the estimate $|\Re \langle \tilde U_E\cdot \nabla w, w\rangle_{L^2_\rho (\R^2)}|\leq C\|\omega\|_{L^2(\R^2)}$. 
Thus we have 
\begin{align}\label{proof.thm.slip.2}
\| (\lambda I - \tilde T_\alpha)^{-1} \tilde g\|_{L^2_\rho (\R^2)} + \frac{1}{\alpha^\frac12} \|\nabla (\lambda I - \tilde T_\alpha)^{-1} \tilde g\|_{L^2(\R^2)}  \leq C \| \tilde g\|_{L^2_\rho(\R^2)}\,,~ \tilde g\in L^2_{\rho,odd}(\R^2)\,.
\end{align}
Here $C$ is independent of $\lambda\in \partial B_\epsilon (\lambda)$ and $\alpha\geq \frac{2}{r_E}$. 
Next we show that, for any $\tilde  g \in L^2_{\rho,odd} (\R^2)$,
\begin{align}\label{proof.thm.slip.3}
\lim_{\alpha\rightarrow \infty} \sup_{\lambda\in \partial B_\epsilon (\lambda_E)} \| (\lambda I - \tilde T_\alpha)^{-1} \tilde g - (\lambda I -\tilde T)^{-1} \tilde g\|_{L^2_\rho (\R^2)} =0\,.
\end{align}
Here $\tilde T \omega = -\nabla \cdot (U_E \omega)$ with $D(\tilde T) = D(\tilde \Lambda_E)$. To prove \eqref{proof.thm.slip.3} we note that the similar energy argument as above for the semigroup $\{e^{\tau \tilde T_\alpha}\}_{\tau\geq 0}$ and $\{e^{\tau \tilde T}\}_{\tau\geq 0}$ in $L^2_{odd}(\R^2)$ and $L^2_{\rho,odd}(\R^2)$ yields
\begin{align}\label{proof.thm.slip.4}
\begin{split}
\| e^{-\lambda \tau} e^{\tau \tilde T_\alpha} \tilde g\|_{L^2(\R^2)} \leq e^{-c_1\tau} \|\tilde g\|_{L^2(\R^2)}\,,\\
\| e^{-\lambda \tau} e^{\tau \tilde T_\alpha} \tilde g \|_{L^2_\rho (\R^2)}\leq C e^{-c_2\tau} \| \tilde g\|_{L^2_\rho (\R^2)}
\end{split}
\end{align}
and 
\begin{align}\label{proof.thm.slip.5}
\begin{split}
\| e^{-\lambda \tau} e^{\tau \tilde T} \tilde g\|_{L^2(\R^2)} \leq e^{-\tau\Re (\lambda)} \|\tilde g\|_{L^2(\R^2)}\,,\\
\| e^{-\lambda \tau} e^{\tau \tilde T} \tilde g \|_{L^2_\rho (\R^2)}\leq C e^{-c'\tau} \| \tilde g\|_{L^2_\rho (\R^2)}
\end{split}
\end{align}
for $\tilde g\in L^2_{\rho,odd}(\R^2)$. Here $C$ and $c_1,c_2,c'$ are positive constants independent of $\lambda\in \partial B_\epsilon (\lambda_E)$ and $\alpha \geq \frac{2}{r_E}$.  
Then we verify from \eqref{proof.thm.slip.4} and \eqref{proof.thm.slip.5} the formula of their resolvents by the Laplace transform 
\begin{align}\label{proof.thm.slip.6}
\begin{split}
(\lambda I - \tilde T_\alpha)^{-1} \tilde g & = \int_0^\infty e^{-\lambda \tau} e^{\tau \tilde T_\alpha} \tilde g \, d\tau\,,\\
(\lambda I - \tilde T)^{-1} \tilde g & = \int_0^\infty e^{-\lambda \tau} e^{\tau \tilde T} \tilde g \, d\tau\,.
\end{split}
\end{align}
Let $\tilde g\in L^2_{\rho,odd}(\R^2)$. For any $n\in \mathbb{N}$ there exists $\tilde g_n \in C_0^\infty (\R^2)\cap L^2_{odd}(\R^2)$ such that $\|\tilde g -\tilde g_n\|_{L^2_{\rho,odd}(\R^2)}<\frac1n$. Then for any $N\geq 1$,
\begin{align*}
& \| (\lambda I - \tilde T_\alpha)^{-1} \tilde g - (\lambda I - \tilde T)^{-1} \tilde g \|_{L^2_{\rho}(\R^2)} \\
& \leq \| (\lambda I - \tilde T_\alpha)^{-1} \tilde g_n - (\lambda I - \tilde T)^{-1} \tilde g_n \|_{L^2_{\rho}(\R^2)}  \\
& \quad + \| (\lambda I - \tilde T_\alpha)^{-1} (\tilde g - \tilde g_n)\|_{L^2_{\rho} (\R^2)} +\|  (\lambda I - \tilde T)^{-1} (\tilde g-\tilde g_n) \|_{L^2_{\rho}(\R^2)} \\
& \leq \int_0^N  \| e^{-\lambda \tau} e^{\tau \tilde T_\alpha} \tilde g_n - e^{-\lambda \tau} e^{\tau \tilde T} \tilde g_n\|_{L^2_\rho (\R^2)} \, d\tau \\
& \quad + \int_N^\infty \big ( \| e^{-\lambda \tau} e^{\tau \tilde T_\alpha} \tilde g_n \|_{L^2_\rho (\R^2)} +\|e^{-\lambda \tau} e^{\tau \tilde T} \tilde g_n\|_{L^2_\rho (\R^2)} \big ) \,d\tau + C \| \tilde g-\tilde g_n \|_{L^2_\rho (\R^2)} \\
& \leq \int_0^N  \| e^{-\lambda \tau} e^{\tau \tilde T_\alpha} \tilde g_n - e^{-\lambda \tau} e^{\tau \tilde T} \tilde g_n\|_{L^2_\rho (\R^2)} \, d\tau \\
& \quad + C(\tilde g)  \int_N^\infty e^{-c\tau} \,d\tau + \frac{C}{n}\,.
\end{align*}
Here the constants $C$ and $c$ are independent of $\lambda\in \partial B_\epsilon (\lambda_E)$ and $\alpha\ge \frac{2}{r_E}$. Since $\tilde g_n$ belongs to $C_0^\infty (\R^2)$, it is easy to see that $\displaystyle \lim_{\alpha\rightarrow \infty} \sup_{0<\tau<N,\lambda\in \partial B_\epsilon (\lambda_E)} \| e^{-\lambda \tau} e^{\tau \tilde T_\alpha} \tilde g_n - e^{-\lambda \tau} e^{\tau \tilde T} \tilde g_n\|_{L^2_\rho (\R^2)} =0$ for each fixed $N\geq 1$. Hence, letting $\alpha \rightarrow \infty$ first and then letting $N,n\rightarrow \infty$, we obtain \eqref{proof.thm.slip.3}. Next we observe that
\begin{align*}
\lambda I - \tilde M_\alpha = (\lambda  I - \tilde T_\alpha ) \big ( I - (\lambda I - \tilde T_\alpha)^{-1} \tilde S \big )\,,
\end{align*}
where $\tilde S[\tilde g] = -\tilde K[\tilde g]\cdot \nabla \Omega_E$ and thus $\tilde S$ is a compact operator in $L^2_{\rho,odd} (\R^2)$. Since $\lambda I - \tilde T_\alpha$ is invertible, the invertibility of $\lambda I - \tilde M_\alpha$ is equivalent with the invertibility of $ I - (\lambda I - \tilde T_\alpha)^{-1} \tilde S$. By the same reason, the invertibility of $\lambda I - \tilde \Lambda_E$ is equivalent with the invertibility of $I-(\lambda I - \tilde T)^{-1} \tilde S$. Since $(\lambda I - \tilde T)^{-1} \tilde S$ is continuous in $\lambda \in \rho_{re} (\tilde T)$ in the operator norm, we have 
\begin{align}\label{proof.thm.slip.7}
\sup_{\lambda\in \partial B_\epsilon(\lambda_E)} \| \big ( I  - (\lambda I - \tilde T)^{-1} \tilde S\big )^{-1}\|_{L^2_{\rho,odd}\rightarrow L^2_{\rho,odd}}<\infty\,.
\end{align}
On the other hand, approximating $\tilde S$ by a finite rank operator in the operator norm, we obtain from \eqref{proof.thm.slip.2} and \eqref{proof.thm.slip.3}, 
\begin{align}\label{proof.thm.slip.8}
\lim_{\alpha\rightarrow \infty} \sup_{\lambda\in \partial B_\epsilon (\lambda_E)} \| (\lambda I - \tilde T_\alpha)^{-1} \tilde S -  (\lambda I - \tilde T)^{-1} \tilde S\|_{L^2_{\rho,odd}\rightarrow L^2_{\rho,odd}} =0\,.
\end{align}
Then \eqref{proof.thm.slip.7} and \eqref{proof.thm.slip.8} imply from the argument of the Neumann series that there exists $\alpha_\epsilon'\geq 1$ such that 
\begin{align}\label{proof.thm.slip.9}
\begin{split}
& \sup_{\lambda\in \partial B_\epsilon(\lambda_E),~\alpha \geq \alpha_\epsilon'} \| \big ( I  - (\lambda I - \tilde T_\alpha)^{-1} \tilde S\big )^{-1}\|_{L^2_{\rho,odd}\rightarrow L^2_{\rho,odd}}<\infty\,,\\
& \lim_{\alpha\rightarrow \infty} \sup_{\lambda\in \partial B_\epsilon (\lambda_E)} \| \big ( I  - (\lambda I - \tilde T_\alpha)^{-1} \tilde S\big )^{-1} - \big ( I  - (\lambda I - \tilde T)^{-1} \tilde S\big )^{-1}\|_{L^2_{\rho,odd}\rightarrow L^2_{\rho,odd}}=0\,.
\end{split}
\end{align}
In particular, $\partial B_\epsilon (\lambda_e) \subset  \rho_{re} (\tilde M_\alpha)$ when $\alpha\geq \alpha_\epsilon'$, and the formua $(\lambda I -\tilde M_\alpha)^{-1} =\big ( I  - (\lambda I - \tilde T_\alpha)^{-1} \tilde S \big )^{-1} (\lambda I - \tilde T_\alpha)^{-1}$ holds, which gives \eqref{est.thm.slip.1} from \eqref{proof.thm.slip.2} and \eqref{proof.thm.slip.9}, and also gives \eqref{est.thm.slip.2} from \eqref{proof.thm.slip.3} and \eqref{proof.thm.slip.9}. The proof is complete.
\end{proof}

Let us recall that $K=\nabla^\bot (-\Delta_D)^{-1}$ is the operator of the Biot-Savart formula in $\R^2_+$. As a direct consequence of Theorem \ref{thm.slip.tilde}, we obtain
\begin{cor}\label{cor.thm.slip.tilde} For any $\epsilon \in (0,\frac{r_E}{2})$ there exists $\alpha_\epsilon'\geq 1$ such that if $\alpha\geq \alpha_\epsilon'$ then $\partial B_\epsilon (\lambda_E) \subset \rho_{re} (M_\alpha)$ and 
\begin{align}\label{est.cor.thm.slip.1}
\begin{split}
& \sup_{\lambda\in \partial B_\epsilon (\lambda_E), \alpha\geq \alpha_\epsilon'}  \big (\| (\lambda I -M_\alpha)^{-1}\|_{L^2_\rho (\R^2_+) \rightarrow L^2_\rho (\R^2_+)} + \| K (\lambda I -M_\alpha)^{-1}\|_{L^2_\rho (\R^2_+) \rightarrow L^2_\sigma (\R^2_+)} \big ) \\
& <\infty\,,
\end{split}
\end{align}
and for any $g\in L^2_{\rho} (\R_+^2)$, 
\begin{align}\label{est.cor.thm.slip.2}
\lim_{\alpha\rightarrow \infty} \sup_{\lambda\in \partial B_\epsilon (\lambda_E)} \|(\lambda I - M_\alpha)^{-1}  g - (\lambda I -  \Lambda_E)^{-1}  g \|_{L^2_{\rho} (\R_+^2)}=0\
\end{align}
\end{cor}

\begin{proof} The assertion 
\begin{align*}
\sup_{\lambda\in \partial B_\epsilon (\lambda_E),\alpha\geq \alpha_\epsilon'} \| (\lambda I -M_\alpha)^{-1}\|_{L^2_\rho \rightarrow L^2_\rho} <\infty
\end{align*}
and \eqref{est.cor.thm.slip.2} directly follow from Theorem \ref{thm.slip.tilde} by the reflection.
As for the estimate of the velocity, it suffices to recall 
\begin{align*}
\| K[h]\|_{L^2(\R^2_+)}\leq C \| \xi_2 h \|_{L^2(\R^2_+)}\,,
\end{align*} 
which is given in Lemma \ref{lem.bs}.
Then the uniform bound of  $K(\lambda I - M_\alpha)^{-1}$ follows from the estimate of the vorticity which is already shown above. The proof is complete. 
\end{proof}

Let $O_{cpt}$ be a smooth bounded domain satisfying ${\rm supp} \,U_E \subset O_{cpt} \subset \overline{O_{cpt}}\subset \R^2_+$. 
For later use we also show the following proposition.

\begin{prop}\label{prop.slip} Let $\epsilon\in (0,\frac{r_E}{2})$, $\alpha \geq \alpha_\epsilon'$, and $\lambda\in \partial B_\epsilon (\lambda_E)$. Assume that $g\in L^2 (\R^2_+)$ satisfies ${\rm supp}\, g\subset \overline{O_{cpt}}$. Then there exists $C>0$ depending only on $\epsilon$ and $O_{cpt}$ such that  
\begin{align}
\| v_1[g]|_{\xi_2=0}\|_{W^{5,2}(\R)}  + \| \xi_1 \partial_1^{l+1} v_1[g]|_{\xi_2=0}\|_{L^2(\R)} \leq C \| g\|_{L^2(\R^2_+)}
\end{align}
for $l=0,1$. Here,
\begin{align*}
v_1[g]=v_1[g;\lambda,\alpha] = \partial_2 (-\Delta_D)^{-1} (\lambda I - M_\alpha)^{-1} g
\end{align*}
is the tangential component of the velocity obtained from the vorticity $(\lambda I - M_\alpha)^{-1} g$ via the Biot-Savart formula. The constant $C$ is independent of $\lambda\in \partial B_\epsilon (\lambda_E)$, $\alpha\geq \alpha_\epsilon'$, and $g$.
\end{prop}

\begin{proof} Let us take $\kappa_0=2^{-6} {\rm dist}\, (\overline{O_{cpt}}, \partial \R^2_+)>0$.
Set $\omega= (\lambda I - M_\alpha)^{-1} g$. We note that $\| \omega\|_{L^2_\rho (\R^2_+)} \leq C\|g\|_{L^2_\rho (\R^2_+)}$ holds by Corollary \ref{cor.thm.slip.tilde}, where $C$ is independent of $\lambda\in \partial B_\epsilon (\lambda_E)$ and $\alpha\geq \alpha_\epsilon'$, and it is also straightforward to see that the integration by part for $\Re \langle (\lambda I - M_\alpha) \omega,\omega\rangle_{L^2_\rho (\R^2_+)} = \Re \langle g,\omega\rangle_{L^2_\rho (\R^2_+)}$ and the $L^2$ estimate of $K[\omega]$ in Lemma \ref{lem.bs} together with the fact that $U_E$ and $\Omega_E$ have a compact support yield the estimate 
\begin{align}\label{proof,prop.slip.1}
\Re (\lambda) \| \omega\|_{L^2_\rho (\R^2_+)} +\frac{1}{\alpha^\frac12} \| \nabla \omega\|_{L^2_\rho (\R^2_+)} & \leq C (\|\omega\|_{L^2_\rho (\R^2_+)} + \| g\|_{L^2_\rho (\R^2_+)}) \nonumber \\
& \leq C \| g\|_{L^2_\rho (\R^2_+)}\,.
\end{align}
Here $C$ is independent of $\lambda\in \partial B_\epsilon (\lambda_E)$ and $\alpha\geq \alpha_\epsilon'$.
Then Lemma \ref{lem.bs} implies 
\begin{align}\label{proof,prop.slip.2}
\begin{split}
& \| v_1[g]|_{\xi_2=0}\|_{W^{5,2}(\R)}  + \sum_{l=0,1} \| \xi_1 \partial_1^{l+1} v_1[g]|_{\xi_2=0}\|_{L^2(\R)} \\
&\leq C\big ( \| \langle \xi\rangle \omega \|_{L^2(\R^2_+)} + \sum_{j=1}^5 \|\partial_1^j \omega\|_{L^2(\R\times (0,\kappa_0))} + \sum_{j=0,1} \| \xi_1\partial_1^{j+1} \omega \|_{L^2(\R\times (0,\kappa_0))}\big )\,.
\end{split}
\end{align}
Since $\|\langle \xi \rangle \omega\|_{L^2(\R^2_+)}\leq C\| \omega\|_{L^2_\rho (\R^2_+)} \leq C\|g\|_{L^2_\rho (\R^2_+)}$, it suffices to show for $\omega =  (\lambda I - M_\alpha)^{-1} g$,
\begin{align}\label{proof,prop.slip.3}
\sum_{j=1}^5\|\partial_1^j \omega\|_{L^2(\R\times (0,\kappa_0))} +  \sum_{l=0,1} \| \xi_1\partial_1^{l+1} \omega \|_{L^2(\R\times (0,\kappa_0))} \leq C\| g\|_{L^2_\rho (\R^2_+)}\,.
\end{align}
Set $\kappa_j = 2^{5-j}\kappa_0$, $1\leq j\leq 5$. Let $\chi_{\kappa_j} = \chi_{\kappa_j} (\xi_2)$ be a smooth cut-off such that $\chi_{\kappa_j}=1$ for $0\leq \xi_2\leq \kappa_j$ and $\chi_{\kappa_j}=0$ for $\xi_2\geq 2\kappa_j$. Since $\omega_j = \partial_1^j \omega$ satisfies $(\lambda I - M_\alpha - \frac{j}{2}) \omega_j = \partial_1^j g$, we see that $\tilde \omega_j =\omega_j \chi_{\kappa_j}$, $j=1,2,3,4,5$, satisfies 
\begin{align}\label{proof,prop.slip.4}
\begin{split}
& \lambda \tilde \omega_j  - \frac{1}{\alpha} (\Delta +\frac{\xi}{2}\cdot \nabla +\frac{2+j}{2}) \tilde \omega_j  = \partial_1 F \,,\\
& F  = -\frac{1}{\alpha} \big ( \chi_{\kappa_j}'' \partial_1^{j-1} \omega +2\chi_{\kappa_j}' \partial_2\partial_1^{j-1}\omega +\frac{\xi_2}{2}\chi_{\kappa_j}' \partial_1^{j-1} \omega \big )\,,
\end{split}
\end{align}
with the boundary condition $\tilde \omega_j|_{\xi_2=0}=0$.
Here we have used that ${\rm supp}\, U_E$ and ${\rm supp}\, g$ are contained in $\overline{O_{cpt}}$ by the assumptions.
Then the integration by part for 
\begin{align*}
\Re \langle \lambda \tilde \omega_j - \frac{1}{\alpha} (\Delta +\frac{\xi}{2}\cdot \nabla +\frac{2+j}{2}) \tilde \omega_j,\tilde \omega_j\rangle_{L^2_\rho (\R^2_+)} = \Re \langle \partial_1 F,\tilde \omega_j\rangle_{L^2_\rho (\R^2_+)}
\end{align*}
and \eqref{est.harmonic} yield
\begin{align}\label{proof,prop.slip.5}
& \Re (\lambda) \| \tilde \omega_j\|_{L^2_\rho(\R^2_+)}  + \frac{1}{\alpha^\frac12} (\| \nabla \tilde \omega_j\|_{L^2_\rho (\R^2_+)} + \| \xi \tilde \omega_j\|_{L^2_\rho (\R^2_+)}) \nonumber\\
& \leq \frac{C}{\alpha^{\frac12}}\|\tilde \omega_j \|_{L^2_\rho (\R^2_+)} + C\alpha^\frac12 \| F\|_{L^2_\rho (\R^2_+)}\nonumber \\ 
\begin{split}
& \leq \frac{C}{\alpha^\frac12} \big ( \|\tilde \omega_j\|_{L^2_\rho (\R^2_+)} +  \| \partial_1^{j-1} \omega \|_{L^2_\rho (\R\times (0,2\kappa_j))} + \| \partial_2 \partial_1^{j-1}\omega\|_{L^2_\rho (\R\times (0,2\kappa_j))}\big )\,.   
\end{split}
\end{align}
This estimate for $j=1$ and \eqref{proof,prop.slip.1} imply 
\begin{align}\label{proof,prop.slip.6}
\| \partial_1 \omega \|_{L^2_\rho(\R\times (0,\kappa_1))}  + \frac{1}{\alpha^\frac12} \| \nabla \partial_1 \omega\|_{L^2_\rho (\R\times (0,\kappa_1))} \leq C\| g\|_{L^2_\rho (\R^2_+)}\,,
\end{align}
and then \eqref{proof,prop.slip.5} for $j=2$ and \eqref{proof,prop.slip.6} with the conditions $2\kappa_2=\kappa_1$ give
\begin{align}\label{proof,prop.slip.7}
\| \partial_1^2 \omega \|_{L^2_\rho(\R\times (0,\kappa_2))}  + \frac{1}{\alpha^\frac12} \| \nabla \partial_1^2 \omega\|_{L^2_\rho (\R\times (0,\kappa_2))} \leq C\| g\|_{L^2_\rho (\R^2_+)}\,.
\end{align}
Similarly, we have from $2\kappa_j=\kappa_{j-1}$ and $\kappa_5=\kappa_0$ that $\|\partial_1^j\omega\|_{L^2_\rho (\R\times (0,\kappa_j))}\leq C\| g\|_{L^2_\rho (\R^2_+)}$ also for $j=3,4,5$. Thus, \eqref{proof,prop.slip.3} holds. The proof is complete.
\end{proof}
 
\section{Linear analysis under noslip boundary condition}\label{sec.noslip}

The aim of this section is to prove Theorem \ref{thm.spectrum}. The next proposition is the key of the proof, in which we construct the solution to the resolvent problem \eqref{eq.resolvent} for a class of given force $f$ with compactly supported vorticity, by making use of the results in Section \ref{sec.slip} and the boundary layer to recover the noslip boundary condition. We recall that $K[\omega]=\nabla^\bot (-\Delta_D)^{-1}\omega$. 

\begin{prop}\label{prop.linear} For any $\epsilon \in (0,\frac{r_E}{2})$ there exists $\alpha_\epsilon''\in [\alpha_\epsilon',\infty)$ such that if $\alpha\geq \alpha_\epsilon''$ then for any $g\in L^2 (\R^2_+)$ with ${\rm supp}\, g\subset \overline{O_{cpt}}$ and $\lambda\in \partial B_\epsilon (\lambda_E)$ there exists a solution $v\in D(\mathbb{L}_\alpha)$ to the resolvent problem $(\lambda I - \mathbb{L}_\alpha) v = K[g]$ satisfying the following condition: $v=v[g;\lambda,\alpha]$ is decomposed as 
\begin{align*}
v [g;\lambda,\alpha]= v_{slip} [g;\lambda,\alpha]+ v_{bc}[g;\lambda,\alpha]\,,
\end{align*}
where 
\begin{align}
v_{slip}[g;\lambda,\alpha] & = K [(\lambda I - M_\alpha )^{-1} g]\,,
\end{align}
while 
\begin{align}
\| v_{bc} [g;\lambda,\alpha] \|_{L^2(\R^2_+)} \leq C \alpha^{-\frac14} \| g\|_{L^2(\R^2_+)}.
\end{align}
Here $C$ depends only on $\epsilon$ and $O_{cpt}$, and is independent of $\lambda\in \partial B_\epsilon (\lambda_E)$, $\alpha\geq \alpha_\epsilon''$, and $g$.
\end{prop}

\begin{proof} Let us construct the solution to the resolvent problem of the form $v=v_{slip} + v_{bc}$, where $v_{slip} =v_{slip}[g]= K [(\lambda I - M_\alpha )^{-1} g]$. Then $v_{bc}$ should obey the system 
\begin{align}\label{eq.v_bc}
\begin{split}
& \lambda v_{bc}  - \frac{1}{\alpha} (\Delta +\frac{\xi}{2}\cdot \nabla +\frac12) v_{bc} + U_E \cdot \nabla v_{bc} + v_{bc}\cdot \nabla U_E + \nabla p_{bc} = 0\,, \quad \xi \in \R_+^2\,,\\
& \nabla\cdot v_{bc}=0,\quad \xi \in \R_+^2\,,\\
& v_{bc,1}|_{\xi_2=0} = h,\quad  v_{bc,2}|_{\xi_2=0} =0\,.
\end{split}
\end{align}
Here $h=-v_{slip,1}[g]|_{\xi_2=0}$. To solve \eqref{eq.v_bc} we intorudce the function space $Z$ for the boundary data $h$ as follows:
\begin{align}\label{proof.prop.linear.1}
\begin{split}
Z & = \{h\in W^{5,2}(\R)~|~\xi_1\partial_1^{l+1}h\in L^2(\R),~l=0,1\}\,,\\
\| h\|_Z & = \| h\|_{W^{5,2}(\R)} + \sum_{l=0,1} \| \xi_1 \partial_1^{l+1} h\|_{L^2(\R)}\,.
\end{split}
\end{align}
We note that Proposition \ref{prop.slip} gives
\begin{align}\label{proof.prop.linear.1'}
\| v_{slip,1}[g]|_{\xi_2=0}\|_{Z} \leq C \| g\|_{L^2 (\R^2_+)}\,,
\end{align}
by recalling the assumption ${\rm supp}\, g\subset \overline{O_{cpt}}\subset \R^2_+$.
For the moment we consider the system \eqref{eq.v_bc} for a given $h\in Z$, rather than the specific one $-v_{slip,1}[g]|_{\xi_2=0}$.
We first correct the boundary data in the form of the boundary layer as follows: 
\begin{align}\label{proof.prop.linear.2}
\begin{split}
J_1[h] (\xi)  & = h (\xi_1) (1-\Xi_2) e^{-\Xi_2}\,, \\
J_2[h] (\xi) & =-  \alpha^{-\frac12} \int_0^{\Xi_2} (1- \eta_2) e^{-\eta_2}\, d\eta_2 \, \partial_1 h(\xi_1)\,,\\
\Xi_2 & =\alpha^{\frac12} \xi_2\,.
\end{split}
\end{align}
It is straightforward to see that
\begin{align*}
\nabla \cdot J[h] =0\,,\quad J_1[h]|_{\xi_2=0} = h\,, \quad J_2[h]|_{\xi_2=0} =0\,,
\end{align*}
and we also have from $\int_0^\infty (1-\eta_2) e^{-\eta_2}\, d\eta_2=0$,
\begin{align}\label{proof.prop.linear.2'}
J_2[h] (\xi) = \alpha^{-\frac12} \int_{\Xi_2}^\infty (1-\eta_2) e^{-\eta_2} \, d\eta_2\, \partial_1 h(\xi_1)\,.
\end{align}
Therefore, $J_2[h]$ decays exponentially as $\Xi_2\rightarrow \infty$. 
Let $O_{cpt}'$ be a smooth bounded domain such that $\overline{O_{cpt}} \subset O_{cpt}' \subset \overline{O_{cpt}'}\subset \R^2_+$. We observe that
\begin{align}\label{proof.prop.linear.3}
\begin{split}
\sum_{k=0,1,2} \alpha^{\frac14-\frac{k}{2}} \| \nabla^k J[h]\|_{L^2(\R^2_+)} + \alpha^{\frac{1}{4}}\| \xi \cdot \nabla J[h] \|_{L^2(\R^2_+)} & \leq C \| h\|_Z\,,\\
\sum_{k=0}^4 \|\nabla^k J[h] \|_{L^2 (O_{cpt}')} & \leq C\alpha^{-1} \| h\|_Z\,, 
\end{split}
\end{align}
which are easily shown from \eqref{proof.prop.linear.2}, \eqref{proof.prop.linear.2'}, and ${\rm dist} (\overline{O_{cpt}'}, \partial\R^2_+)>0$.
Next we set
\begin{align}\label{proof.prop.linear.4}
R [h] =  (\lambda I - \mathbb{H}_\alpha )^{-1} \Big ( -\lambda  +  \frac{1}{\alpha} (\mathbb{P}_\sigma \Delta +\frac{\xi}{2}\cdot \nabla +\frac12) \Big )J[h]\,,
\end{align}
where $\mathbb{H}_\alpha = \frac{1}{\alpha} ( \mathbb{A} + \frac{\xi}{2}\cdot \nabla +\frac12)$, by noticing $J[h], \frac{\xi}{2}\cdot \nabla J[h]\in L^2_\sigma (\R^2_+)$. 
Then $R[h]$ belongs to $D(\mathbb{L}_\alpha)$, and Lemma \ref{lem.stokes} yields 
\begin{align}\label{proof.prop.linear.5}
\| R[h]\|_{L^2(\R^2_+)} & \leq C \| \big ( -\lambda +  \frac{1}{\alpha} (\mathbb{P}_\sigma \Delta +\frac{\xi}{2}\cdot \nabla +\frac12) \big ) J[h] \|_{L^2(\R^2_+)} \nonumber \\
& \leq C\alpha^{-\frac14} \| h\|_Z\,,
\end{align}
where \eqref{proof.prop.linear.3} is used in the last line.
We also have from Lemma \ref{lem.stokes} that 
\begin{align}\label{proof.prop.linear.5'}
& \sum_{k=0,1} \| \nabla^k {\rm rot}\, R[h]\|_{L^2(O_{cpt})} \nonumber \\
& \leq C \| \big ( -\lambda +  \frac{1}{\alpha} (\mathbb{P}_\sigma \Delta +\frac{\xi}{2}\cdot \nabla +\frac12) \big ) J[h] \|_{L^2(\R^2_+)} \nonumber \\
& \quad + C\sum_{k=0,1}\| \nabla^k {\rm rot}\, \big ( -\lambda +  \frac{1}{\alpha} (\mathbb{P}_\sigma \Delta +\frac{\xi}{2}\cdot \nabla +\frac12) \big ) J[h] \|_{L^2(O_{cpt}')}\nonumber \\
& \leq C \alpha^{-\frac14} \| h\|_Z\,.
\end{align}
Here we have also used ${\rm rot} \,\mathbb{P}_\sigma = {\rm rot}$ and \eqref{proof.prop.linear.3} in the last line.
Note that $w = J[h] + R[h]$ satisfies 
\begin{align*}
\begin{split}
& \lambda w - \frac{1}{\alpha} (\Delta +\frac{\xi}{2}\cdot \nabla +\frac12) w + U_E \cdot \nabla w + w \cdot \nabla U_E + \nabla q = U_E \cdot \nabla w + w \cdot \nabla U_E \,,\\
& \nabla\cdot w=0\,,\\
& w_{1}|_{\xi_2=0} = h,\quad  w_{2}|_{\xi_2=0} =0\,.
\end{split}
\end{align*}
Then we set 
\begin{align}\label{proof.prop.linear.7}
v_{slip}^{(2)} [h] = - K[(\lambda I - M_\alpha)^{-1} {\rm rot}\, \big ( U_E \cdot \nabla w + w \cdot \nabla U_E\big )]\,,
\end{align}
which implies that $v_{app} [h]= J[h]+R[h] + v_{slip}^{(2)}[h]$ satisfies 
\begin{align}\label{proof.prop.linear.8}
\begin{split}
& \lambda v_{app} - \frac{1}{\alpha} (\Delta +\frac{\xi}{2}\cdot \nabla +\frac12) v_{app} + U_E \cdot \nabla v_{app} + v_{app} \cdot \nabla U_E + \nabla p_{app}= 0 \,,\\
& \nabla\cdot v_{app}=0\,,\\
& v_{app,1}|_{\xi_2=0} = h + \Psi[h],\quad  v_{app,2}|_{\xi_2=0} =0\,,
\end{split}
\end{align}
Here $\Psi[h] = v^{(2)}_{slip,1}[h]|_{\xi_2=0}$. 
For the estimate of $v_{slip}^{(2)}$, the key point is that, since ${\rm supp}\, U_E\subset O_{cpt}\subset \overline{O_{cpt}}\subset \R^2_+$, the term 
\begin{align*}
{\rm rot}\, \big ( U_E \cdot \nabla w + w \cdot \nabla U_E\big ) = U_E \cdot \nabla {\rm rot}\, (J+R)[h] + (J+R)[h] \cdot \nabla \Omega_E
\end{align*}
is small and has the order $O(\alpha^{-\frac14})$ for its $L^2$ norm; indeed, by using \eqref{proof.prop.linear.3}, \eqref{proof.prop.linear.5}, and \eqref{proof.prop.linear.5'}, we see that there exists $C>0$ independent of $\lambda \in \partial B_\epsilon (\lambda_E)$ and $\alpha$ such that 
\begin{align}\label{proof.prop.linear.9}
\| {\rm rot}\, \big ( U_E \cdot \nabla w + w \cdot \nabla U_E\big ) \|_{L^2(\R^2_+)} \leq C\alpha^{-\frac14} \| h\|_Z\,.
\end{align}
Hence we have from Corollary \ref{cor.thm.slip.tilde} together with the condition that $U_E$ has a compact support,
\begin{align}\label{proof.prop.linear.10}
\| v_{slip}^{(2)} [h]\|_{L^2(\R^2_+)} + \| {\rm rot}\, v_{slip}^{(2)}[h]\|_{L^2_\rho (\R^2_+)} & \leq C\| {\rm rot}\, \big ( U_E \cdot \nabla w + w \cdot \nabla U_E\big ) \|_{L^2(\R^2_+)}  \nonumber \\
& \leq C \alpha^{-\frac14} \| h\|_Z\,,
\end{align}
while Proposition \ref{prop.slip} implies
\begin{align}\label{proof.prop.linear.11}
\| \Psi[h]\|_Z \leq C\| {\rm rot}\, \big ( U_E \cdot \nabla w + w \cdot \nabla U_E\big ) \|_{L^2(\R^2_+)} \leq C\alpha^{-\frac14} \| h\|_Z\,.
\end{align}
Here $C$ is independent of $\lambda\in \partial B_\epsilon (\lambda_E)$ and $\alpha\in [\alpha_\epsilon',\infty)$. 
Therefore, the linear map $I+\Psi$ is invertible in $Z$ for large enough $\alpha$ (independently of $\lambda$), and its inverse is given by the Neumann series. Note that we can take $\alpha_\epsilon''$ large enough so that $\|(I+\Psi)^{-1}\|_{Z\rightarrow Z}\leq 2$. The solution $v_{bc}=v_{bc}[h]$ of \eqref{eq.v_bc} is then given by the formula 
\begin{align*}
v_{bc}[h] & = v_{app} [ (I+\Psi)^{-1} h]  \\
& = J[ (I+\Psi)^{-1} h] + R[ (I+\Psi)^{-1} h] + v_{slip}^{(2)} [ (I+\Psi)^{-1} h],
\end{align*}
which satisfies, when $h=-v_{slip,1}[g]|_{\xi_2=0}$, by recalling \eqref{proof.prop.linear.1'},
\begin{align}\label{proof.prop.linear.12}
\| v_{bc}[h]\|_{L^2 (\R^2_+)}\leq C\alpha^{-\frac14} \|h\|_Z \leq C \alpha^{-\frac14} \| g\|_{L^2(\R^2_+)}\,.
\end{align}
From the construction,  it is easy to check that our solution $v=v_{slip}[g] + v_{bc}[-v_{slip,1}[g]]$ belongs to $D(\mathbb{L}_\alpha)$. The proof is complete.
\end{proof}

\begin{proof}[Proof of Theorem \ref{thm.spectrum}.] Suppose that the assertion of Theorem \ref{thm.spectrum} does not hold. Then there exists $\epsilon_0\in (0, \frac{r_E}{2})$ such that there exists a sequence $\{\alpha_n\}$, $\lim_{n\rightarrow \infty} \alpha_n =\infty$, such that $B_{\epsilon_0} (\lambda_E)$ does not contain isolated eigenvalues of $\mathbb{L}_{\alpha_n}$. Since possible spectrum of $\mathbb{L}_{\alpha_n}$ in $B_{\epsilon_0}(\lambda_E)$ consists only of isolated eigenvalues, we conclude that $B_{\epsilon_0} (\lambda_E) \subset \rho_{re} (\mathbb{L}_{\alpha_n})$. Then we have 
\begin{align}\label{proof.thm.spectrum.1}
\frac{1}{2\pi i} \int_{\partial B_{\frac{\epsilon_0}{2}} (\lambda_E)} (\lambda I - \mathbb{L}_{\alpha_n})^{-1} d\lambda =0, \quad n\in \mathbb{N}\,.
\end{align}
Take any $g\in L^2(\R^2_+)$  with ${\rm supp}\, g\subset \overline{O_{cpt}}$, where the bounded domain $O_{cpt}$ is taken so that ${\rm supp}\, U_E \subset O_{cpt} \subset \overline{O_{cpt}} \subset \R^2_+$. Then Proposition \ref{prop.linear} implies that, for any $\lambda\in \partial B_{\frac{\epsilon_0}{2}} (\lambda_E)$,
\begin{align}\label{proof.thm.spectrum.2}
(\lambda I - \mathbb{L}_{\alpha_n})^{-1} K[g] = K[(\lambda I - M_{\alpha_n} )^{-1} g] + v_{bc}[g;\lambda,\alpha_n]\,.
\end{align}
Since $\| K[h]\|_{L^2(\R^2_+)}\leq C\| \xi_2 h\|_{L^2(\R^2_+)} \leq C \|h\|_{L^2_\rho (\R^2_+)}$, in virtue of Corollary \ref{cor.thm.slip.tilde}, the first term in the right-hand side of \eqref{proof.thm.spectrum.2} converges to $K[ (\lambda I -\Lambda_E)^{-1} g]$ in $L^2_\sigma (\R^2_+)$ as $n\rightarrow \infty$ uniformly in $\lambda\in \partial B_{\frac{\epsilon_0}{2}} (\lambda_E)$. On the other hand, Proposition \ref{prop.linear} implies $v_{bc}[g;\lambda,\alpha_n]$ converges to $0$ in $L^2_\sigma (\R^2_+)$ uniformly in $\lambda\in \partial B_{\frac{\epsilon_0}{2}}(\lambda_E)$. 
Hence, we have
\begin{align*}
0 & = \frac{1}{2\pi i} \int_{\partial B_{\frac{\epsilon_0}{2}} (\lambda_E)} \Big ( K[(\lambda I - M_{\alpha_n} )^{-1} g ] + v_{bc}[g;\lambda,\alpha_n] \Big ) d\lambda \nonumber \\
& \rightarrow \frac{1}{2\pi i} \int_{\partial B_{\frac{\epsilon_0}{2}} (\lambda_E)} K[(\lambda I - \Lambda_E )^{-1} g]  \, d\lambda ~~{\rm in}~L^2_\sigma (\R^2_+) \quad (n\rightarrow \infty)\,.
\end{align*}
Thus, $ \frac{1}{2\pi i} \int_{\partial B_{\frac{\epsilon_0}{2}} (\lambda_E)}  K[ (\lambda I - \Lambda_E )^{-1} g]  \, d\lambda =0$ holds, and then, by taking ${\rm rot}$, we obtain $ \frac{1}{2\pi i} \int_{\partial B_{\frac{\epsilon_0}{2}} (\lambda_E)}  (\lambda I - \Lambda_E )^{-1} g  \, d\lambda =0$ as well. Since $g\in L^2(\R^2_+)$ with ${\rm supp}\, g\subset \overline{O_{cpt}}$ is arbitrary, this implies 
\begin{align*}
\frac{1}{2\pi i} \int_{\partial B_{\frac{\epsilon_0}{2}} (\lambda_E)}  (\lambda I - \Lambda_E|_{O_{cpt}} )^{-1}   \, d\lambda =0\,,
\end{align*}
which is a contradiction since $\lambda_E$ is an isolated eigenvalue of $\Lambda_E|_{O_{cpt}}$ and $B_{\epsilon_0} (\lambda_E) \setminus \{\lambda_E\}\subset \rho_{re} (\Lambda_E|_{O_{cpt}})$; see Remark \ref{rem.assume}.  
The proof of Theorem \ref{thm.spectrum} is complete.
\end{proof}

\

We next show the estimates of the semigroup $\{e^{\tau \alpha \mathbb{L}_\alpha}\}_{\tau \geq 0}$ generated by the operator $\alpha \mathbb{L}_\alpha = \mathbb{H} - \alpha \mathbb{P}_\sigma (U_E \cdot \nabla f+ f\cdot \nabla U_E)$ in $L^2_\sigma (\R^2_+)$.

\begin{prop}\label{prop.semigroup} Let $\alpha_E =\alpha_\epsilon$ be the number in Theorem \ref{thm.spectrum} for some $\epsilon\in (0,\frac{r_E}{2})$, and let $\alpha \geq \alpha_E$. Then there exists an unstable eigenvalue $\lambda_{\alpha,max}$ of $\mathbb{L_\alpha}$ in $L_\sigma^2(\R^2_+)$ such that 
\begin{align}\label{est.prop.semigroup.1}
\Re (\lambda_{\alpha,max}) = \sup_{\lambda \in \sigma (\mathbb{L}_\alpha)} \Re (\lambda) >0\,.
\end{align}  
Moreover, for any $\delta>0$ there exists $C_\delta>0$ such that for any $f\in L^2_\sigma (\R^2_+)$ and $\tau>0$, 
\begin{align}\label{est.prop.semigroup.2}
\begin{split}
& \| e^{\tau \alpha \mathbb{L}_\alpha} f\|_{L^2(\R^2_+)} + a(\tau)^\frac12 \| \nabla e^{\tau \alpha \mathbb{L}_\alpha} f\|_{L^2(\R^2_+)} + a(\tau)^\frac12 \| e^{\tau \alpha \mathbb{L}_\alpha} f\|_{L^\infty(\R^2_+)} \\
& \leq C_\delta e^{(\alpha\Re (\lambda_{\alpha,max}) +\delta) \tau} \| f\|_{L^2(\R^2_+)}\,.
\end{split}
\end{align}
Here $a(\tau) = 1-e^{-\tau}$.
\end{prop}

\begin{rem}{\rm As for the maximal eigenvalue $\lambda_{\alpha,max}$ of $\mathbb{L}_\alpha$ above, we can show that 
\begin{align*}
0<\inf_{\alpha \geq \alpha_E} \Re (\lambda_{\alpha,max}) \leq \sup_{\alpha \geq \alpha_E} \Re (\lambda_{\alpha,max}) <\infty\,,
\end{align*}
although we do not need this uniform bound in this paper. 
}
\end{rem}

\begin{proof} We first note that \eqref{est.prop.semigroup.2} for $0<\tau \leq \tau_0$ for small enough $\tau_0>0$ is obtained by solving the integral equation for $v(\tau) = e^{\tau \alpha \mathbb{L}_\alpha} f$, $f\in L^2_\sigma (\R^2_+)$, as follows:
\begin{align}\label{proof.prop.semigroup.1}
\begin{split}
v(\tau) & = e^{\tau \mathbb{H}} f  + \alpha \Big (\Psi [U_E, v] (\tau)  + \Psi [v,U_E](\tau) \Big )\,,\\
\Psi[f,g] (\tau) & = - \int_0^\tau e^{(\tau-s)\mathbb{H}} \mathbb{P}_\sigma \nabla \cdot (f(s) \otimes g (s)) \, ds\,.
\end{split}
\end{align}
Indeed, Lemma \ref{lem.stokes.semigroup} yields 
\begin{align*}
\| \Psi[v,\tilde v] (\tau)\|_{L^2(\R_+^2)} & \leq C\int_0^\tau \| v(s)\|_{L^\infty(\R_+^2)} \| \nabla \tilde v(s)\|_{L^2(\R^2_+)}\, ds
\end{align*}
and 
\begin{align*}
\begin{split}
& \| \nabla \Psi[v,\tilde v] (\tau)\|_{L^2(\R_+^2)} + \| \Psi[v,\tilde v] (\tau)\|_{L^\infty(\R^2_+)}  \\
& \leq C\int_0^\tau  \frac{1}{a(\tau-s)^\frac12} \| v(s)\|_{L^\infty(\R_+^2)} \| \nabla \tilde v(s)\|_{L^2(\R^2_+)}\, ds\,,
\end{split}
\end{align*}
which implies 
\begin{align}\label{proof.prop.semigroup.4}
\| \Psi [U_E, v] \|_{Y_{\tau_0}} + \| \Psi [v,U_E] \|_{Y_{\tau_0}} \leq C \tau_0^\frac12 \| v\|_{Y_{\tau_0}}\,,
\end{align}
where 
\begin{align*}
\| v\|_{Y_{\tau_0}} = \sup_{0<\tau\leq \tau_0} \Big ( \| v(\tau) \|_{L^2(\R^2_+)} + a(\tau)^\frac12 \| \nabla v(\tau)\|_{L^2(\R^2_+)} + a(\tau)^\frac12 \| v(\tau)\|_{L^\infty (\R^2_+)}\Big )\,.
\end{align*}
Note that Lemma \ref{lem.stokes.semigroup} implies that
\begin{align}\label{proof.prop.semigroup.6}
\| e^{\cdot \mathbb{H}} f\|_{Y_{\tau_0}} \leq C \|f\|_{L^2(\R^2_+)},\quad f\in L^2_\sigma (\R^2_+)\,.
\end{align}
Here $C$ is independent of $\tau_0$ and $f\in L^2_\sigma (\R^2_+)$.
Hence, by setting the closed ball $Y_{\tau_0, R_0} = \{ v\in C([0,\tau_0]; L^2_\sigma (\R^2_+))~|~a(\tau)^\frac12 \nabla v (\tau)\in L^\infty (0,\tau_0; L^2(\R_+^2))\,, a(\tau)^\frac12 v(\tau)\in L^\infty (0,\tau_0; L^\infty (\R^2_+))\,, \| v\|_{Y_{\tau_0}}\leq R_0 \}$ with $R_0 = 2\| e^{\cdot \mathbb{H}} f\|_{Y_{\tau_0}}$, we see that the map 
\begin{align}\label{proof.prop.semigroup.7}
Y_{\tau_0,R_0}\ni v \mapsto e^{\tau \mathbb{H}} f +\alpha  \Big (\Psi [U_E, v] (\tau)  + \Psi [v,U_E](\tau) \Big ) \in Y_{\tau_0,R_0}
\end{align}
is well-defined and is a contraction if $\tau_0$ is small enough, where such $\tau_0$ is taken uniformly in $f\in L^2_\sigma (\R^2_+)$. Hence, the unique fixed point in $Y_{\tau_0,R_0}$ defines the strongly continuous semigroup $e^{\tau \alpha \mathbb{L}_\alpha}$ in $L^2_\sigma (\R^2_+)$ for $0<\tau\leq \tau_0$, and by repeating this argument for $f$ replaced by $e^{\alpha \tau_0\mathbb{L}_\alpha} f$, we obtain the whole semigroup $e^{\tau \alpha \mathbb{L}_\alpha}$, $\tau\geq 0$. 
To complete the proof of \eqref{est.prop.semigroup.1} and \eqref{est.prop.semigroup.2}, by the semigroup property, it then suffices to show that the growth bound of $\{e^{\tau \alpha \mathbb{L}_\alpha}\}_{\tau\geq 0}$ is given by $\alpha \Re (\lambda_{\alpha,max})$, i.e., it suffices to show 
\begin{align}\label{proof.prop.semigroup.8}
\begin{split}
& \alpha \Re (\lambda_{\alpha,max}) \\
& = \inf \{ \gamma \in \R~|~ {\rm there ~exists~} C_\gamma>0~{\rm such~that}\\
& \qquad \qquad \| e^{\tau \alpha \mathbb{L}_\alpha} f\|_{L^2(\R^2_+)}\leq C_\gamma e^{\gamma \tau} \|f\|_{L^2(\R^2_+)}~{\rm for~any}~\tau\geq 0~{\rm and}~f\in L^2_\sigma (\R^2_+)\}\,.
\end{split}
\end{align}
To show this we first observe that $e^{\tau \alpha \mathbb{L}_\alpha}$ is written as $e^{\tau \mathbb{H}} + \alpha S(\tau)$, see \eqref{proof.prop.semigroup.1}, and $S(\tau)$ is a compact operator  in $L^2_\sigma (\R^2_+)$. Indeed, we have 
\begin{align}\label{proof.prop.semigroup.8'}
& \| S(\tau) f\|_{L^2(\R^2_+)} + \| \nabla S(\tau) f\|_{L^2(\R^2_+)} \nonumber \\
& \leq  C\int_0^\tau \frac{1}{a(\tau-s)^\frac12} (\| \nabla e^{s\alpha \mathbb{L}_\alpha} f\|_{L^2(\R^2_+)} + \| e^{s\alpha \mathbb{L}_\alpha} f\|_{L^\infty(\R^2_+)}) \, ds \nonumber \\
& \leq C_\tau  \int_0^\tau \frac{1}{a(\tau-s)^\frac12 a(s)^\frac12}  \|f\|_{L^2(\R^2_+)} \, ds\nonumber \\
& \leq C_\tau \| f\|_{L^2(\R^2_+)}\,,
\end{align}
and,  by setting $v(s) = e^{s\alpha \mathbb{L}_\alpha} f$, we also have from Lemma \ref{lem.stokes.semigroup} that
\begin{align}\label{proof.prop.semigroup.8''}
\| \langle \xi\rangle^\frac14 S(\tau) \|_{L^2(\R^2_+)} & \leq \int_0^\tau \| \langle e^\frac{\tau-s}{2} \xi\rangle^\frac14  e^{(\tau-s)\mathbb{H}} \mathbb{P}_\sigma \nabla \cdot \Big (U_E\otimes v (s) + v(s) \otimes U_E)\Big ) \|_{L^2(\R^2_+)}\, ds\nonumber \\
& \leq C e^\frac{\tau}{4} \int_0^\tau \| \langle \xi\rangle^\frac14 \mathbb{P}_\sigma \nabla \cdot \Big (U_E\otimes v (s) + v(s) \otimes U_E)\Big ) \|_{L^2(\R^2_+)}\, ds\nonumber \\
& \leq C_\tau \int_0^\tau  \| \langle \xi\rangle^\frac14 \nabla \cdot \Big (U_E\otimes v (s) + v(s) \otimes U_E)\Big ) \|_{L^2(\R^2_+)}\, ds\nonumber \\
& \leq C_\tau \int_0^\tau (\| \nabla v(s)\|_{L^2(\R^2_+)} + \| v(s) \|_{L^2(\R^2_+)}) \, ds\nonumber \\
& \leq C_\tau \| f\|_{L^2(\R^2_+)}\,.
\end{align}
Here the constant $C_\tau$ depends only on $\tau$. In the third line of \eqref{proof.prop.semigroup.8''} we have used the fact that the Helmholtz projection $\mathbb{P}_\sigma$ is bounded also in the weighted $L^2$ space, see \cite[Theorem 2.4]{KK}. The estimates \eqref{proof.prop.semigroup.8'} and \eqref{proof.prop.semigroup.8''} imply that $S(\tau)$ is a compact operator in $L^2_\sigma (\R^2_+)$.
Let $r_{ess}(e^{\tau \alpha \mathbb{L}_\alpha})$ be the radius of the essential spectrum of $e^{\tau \alpha \mathbb{L}_\alpha}$. Since $S(\tau)$ is compact, we have $r_{ess} (e^{\tau \alpha \mathbb{L}_\alpha} )= r_{ess} (e^{\tau \mathbb{H}}) \leq 1$. This implies from \cite[Corollary 2.11]{EN} that the set of the unstable spectrum $\sigma_{unst} (\alpha \mathbb{L}_\alpha) := \{\lambda\in \sigma (\alpha \mathbb{L}_\alpha)~|~\Re (\lambda)>0\}$ consists only of isolated eigenvalues with finite algebraic multiplicities and the set $\{\lambda\in \sigma (\alpha \mathbb{L}_\alpha)~|~\Re (\lambda)>\delta\}$ is finite for any $\delta>0$, and as long as $\sigma_{unst} (\alpha \mathbb{L}_\alpha)$ is not empty (this condition is ensured by Theorem \ref{thm.spectrum}), the spectral bound $s(\alpha \mathbb{L}_\alpha):=\sup \,\{ \Re (\lambda) ~|~\lambda\in \sigma (\alpha \mathbb{L}_\alpha)\} = \max\, \{ \Re (\lambda)~|~\lambda \in \sigma_{unst} (\alpha \mathbb{L}_\alpha)\}$ coincides with the growth bound of the semigroup $\{e^{\tau \alpha \mathbb{L}_\alpha}\}_{\tau\geq 0}$, which verifies \eqref{proof.prop.semigroup.8}. The proof is complete.
\end{proof}

\section{Construction of solutions to the nonlinear problem}\label{sec.nonlinear}

In this section we prove Theorem \ref{thm.nonlinear}. The argument is almost parallel to that in \cite[Section 4]{ABC1}. Since the self-similar velocity $\alpha u_E = \frac{\alpha}{\sqrt{t}} U_E (\frac{x}{\sqrt{t}})$ is clearly a mild solution of \eqref{eq.ns} with the force $F=F_\alpha = \alpha (\partial_t u_E-\Delta u_E)$, it suffices to construct another solution different from $\alpha u_E$. Let $\alpha \geq \alpha_E$, where $\alpha_E$ is the number in Proposition \ref{prop.semigroup}.  
To simplify the notation, we write $\lambda_\alpha$ for $\lambda_{\alpha,max}$ of Proposition \ref{prop.semigroup}. Then we will construct the solution $u$ of \eqref{eq.ns} in the form 
\begin{align}\label{proof.thm.nonlinear.1}
u(t) = \frac{\alpha}{\sqrt{t}} U_E (\frac{x}{\sqrt{t}}) + \frac{1}{\sqrt{t}} v_{unst} (\log t,\frac{x}{\sqrt{t}}) + \frac{1}{\sqrt{t}} w (\log t,\frac{x}{\sqrt{t}})\,,
\end{align}
where 
\begin{align}\label{proof.thm.nonlinear.2}
v_{unst} (\tau,\xi) = \Re \big ( e^{\alpha \lambda_\alpha \tau}  V_{unst} (\xi)\big )\,, \quad \lambda_{\alpha} = \lambda_{\alpha,max}\,,
\end{align}
and $V_{unst}\in D(\mathbb{L}_\alpha)$ is an eigenfunction for the maximal unstable eigenvalue $\alpha \lambda_\alpha$ of the operator $\alpha \mathbb{L}_\alpha$. We may normalize the norm of the profile $V_{unst}$ so that $\| V_{unst} \|_{W^{2,2}(\R^2_+)} =1$. Thus, our aim is to construct the remainder $w(\tau,\xi)$ by solving the integral equations in the self-similar variables
\begin{align}\label{proof.thm.nonlinear.3}
\begin{split}
w(\tau) & = \Phi [v_{unst},v_{unst}](\tau) + \Phi [v_{unst},w](\tau) + \Phi [w,v_{unst}](\tau) + \Phi [w,w](\tau)\,,\\
\Phi[f,h](\tau) & = -\int_{-\infty}^\tau e^{(\tau-s)\alpha \mathbb{L}_\alpha} \mathbb{P}_\sigma \nabla \cdot (f\otimes g)(s)\, ds\,.
\end{split}
\end{align}
As usual, this can be done by the standard Banach fixed point theorem in the closed ball 
\begin{align}\label{proof.thm.nonlinear.4}
\begin{split}
X_{R} & = \big \{f\in C((-\infty, -R]; L^2_\sigma (\R^2_+)\cap W^{1,2}_0(\R^2_+)^2) ~\big |~\| f\|_{X_R} \leq 1 \big \}\,, \\
\| f\|_{X_R} & = \sup_{\tau\in (-\infty,-R]} e^{-\frac32\alpha \Re(\lambda_\alpha) \tau} \big (\|f(\tau)\|_{W^{1,2} (\R^2_+)} + \|f(\tau)\|_{L^\infty (\R^2_+)} \big )\,.
\end{split}
\end{align}
Here $R\geq 1$ will be taken large enough. Fix $\delta\in (0,1)$ small enough so that $\delta<\frac12 \alpha \Re (\lambda_\alpha)$ holds. Then, Proposition \ref{prop.semigroup} implies 
\begin{align}\label{proof.thm.nonlinear.5}
& \| \Phi [v_{unst},v_{unst}] \|_{X_R} \nonumber \\
& \leq C \sup_{\tau \in (-\infty,-R]} e^{-\frac32\alpha \Re(\lambda_\alpha) \tau} \int_{-\infty}^\tau \frac{1}{a(\tau-s)^\frac12} e^{(\alpha \Re (\lambda_\alpha) +\delta)(\tau-s)} e^{2\alpha \Re (\lambda_\alpha)s} \, ds \nonumber \\
& \qquad \times \|V_{unst}\|_{W^{2,2}(\R^2_+)}^2 \nonumber \\
& \leq C e^{-\frac12\alpha \Re (\lambda_\alpha) R}\,,
\end{align}
\begin{align}\label{proof.thm.nonlinear.5'}
& \| \Phi [v_{unst},f] \|_{X_R} +  \| \Phi [f,v_{unst}] \|_{X_R} \nonumber \\
& \leq C \sup_{\tau \in (-\infty,-R]} e^{-\frac32\alpha \Re(\lambda_\alpha) \tau} \int_{-\infty}^\tau \frac{1}{a(\tau-s)^\frac12} e^{(\alpha \Re (\lambda_\alpha) +\delta)(\tau-s)} e^{\frac52\alpha \Re (\lambda_\alpha)s} \, ds\nonumber \\
& \qquad \times  \|V_{unst}\|_{W^{2,2}(\R^2_+)} \|f\|_{X_R} \nonumber \\
& \leq C e^{-\alpha \Re (\lambda_\alpha) R} \| f\|_{X_R}\,,
\end{align}
and 
\begin{align}\label{proof.thm.nonlinear.6}
& \| \Phi [f,g] \|_{X_R} \nonumber \\
& \leq C \sup_{\tau \in (-\infty,-R]} e^{-\frac32\alpha \Re(\lambda_\alpha) \tau} \int_{-\infty}^\tau \frac{1}{a(\tau-s)^\frac12} e^{(\alpha \Re \lambda_\alpha +\delta)(\tau-s)} e^{3\alpha \Re (\lambda_\alpha)s} \, ds \| f\|_{X_R} \| g\|_{X_R} \nonumber \\
& \leq C e^{-\frac32\alpha \Re (\lambda_\alpha) R}\|f\|_{X_R} \| g\|_{X_R}
\end{align}
for $f,g\in X_R$. Therefore, by taking $R$ large enough, the map $X_R\ni w \mapsto   \Phi [v_{unst},v_{unst}] +  \Phi [v_{unst},w] +  \Phi [w,v_{unst}] + \Phi [w,w]\in X_R$ is a contraction, and hence, there exists a unqiue solution $w\in X_R$ of \eqref{proof.thm.nonlinear.3}.
We note that $v_{unst} (\tau)+w(\tau)$ is not identically $0$, for $w(\tau)$ decays faster than $v_{unst}(\tau)$ as $\tau\rightarrow -\infty$. We also note that, once we construct $u$ of the form \eqref{proof.thm.nonlinear.1} for $(0,t_0]$ with some $t_0>0$, then it is extended as a global (in time) solution by the standard theory of the two-dimensional Navier-Stokes equations in the energy class.  In particular, $u$ is a mild solution of \eqref{eq.ns} in the sense of Definition \ref{def.mild}. Indeed, the properties (i) and (iii) are satisfied by the construction, and as for (ii), we see that $\frac{\alpha}{\sqrt{t}} U_E (\frac{x}{\sqrt{t}}) \rightarrow 0$ as $t\rightarrow 0$ weakly in $L^2_\sigma(\R^2_+)$, while $\frac{1}{\sqrt{t}} v_{unst} (\log t,\frac{x}{\sqrt{t}}) + \frac{1}{\sqrt{t}} w (\log t,\frac{x}{\sqrt{t}})$ converges to $0$ as $t\rightarrow 0$ strongly in $L^2_\sigma (\R^2_+)$ by the construction. The proof of Theorem \ref{thm.nonlinear} is complete.

\section{Example of unstable stationary Euler flows in $\R^2_+$}\label{sec.euler}

In this section, we construct an example of unstable stationary Euler flows in $\R^2_+$.
For this purpose, we consider a stationary, skew-symmetric vortex on the $\xi_1$-axis in $\mathbb{R}^2$. 
Before constructing the example, we recall the instability result by Albritton, Bru{\'e}, and Colombo \cite{ABC1} and check that the linear operator remains unstable in $L^2$ spaces. 
We first recall the instability result by \cite{ABC1}. 
For integer $m \ge 2$, let $L^2_{m}(\mathbb{R}^2)$ be the space of the m-fold rotationally symmetric fuctions:
\[
\begin{aligned}
    L^2_m (\mathbb{R}^2)
    :=
    \{
        f \in L^2 (\mathbb{R}^2)| f(R_{\frac{2\pi}{m}}\xi) = f(\xi)   \quad {\rm for \ a.e. \, }  \xi \in \mathbb{R}^2
        \},
\end{aligned}
\]
where $R_{\frac{2\pi}{m}} : \mathbb{R}^2 \to \mathbb{R}^2$ is counterclockwise rotation by $2\pi/m$.
Let us consider a smooth, divergence-free vector field $\bar{u} (\xi) = \zeta (|\xi|) \xi^{\bot}$ such that $\zeta$ decays at infinity and $\xi^{\bot}=(-\xi_2, \xi_1)$. 
We set $\bar{\omega} =\nabla^{\bot} \cdot \bar{u}$ which is a radial function. 
We then define the linear operator $\tilde \Lambda_{0}$ in $L^2_m(\mathbb{R}^2)$ as follows:
\[
\begin{aligned}
    \tilde \Lambda_{0}: 
    D(\tilde \Lambda_{0}) \subset 
    L^2_m (\mathbb{R}^2) \to L^2_m (\mathbb{R}^2), 
    &\quad
    D(\tilde \Lambda_{0})
    := 
    \{
        \omega \in L^2_m (\mathbb{R}^2) |\nabla \cdot (\bar{u} \omega) \in L^2_m (\mathbb{R}^2)
    \},
    \\
    \tilde \Lambda_{0}\omega 
    &= 
    - \nabla \cdot (\bar{u} \omega)
    - (\tilde K [\omega] \cdot \nabla) \bar{\omega}.
\end{aligned}
\]
Here, we recall that $\tilde K [\omega] = \nabla^{\bot} (-\Delta)^{-1} \omega$ is the velocity in $\mathbb{R}^2$ obtained from the vorticity $\omega$ via the Biot--Savart formula.
Vishik \cite{V2} shows the instability of the linear operator $\tilde \Lambda_{0}$ with an appropriate vorticity profile $\bar{\omega}$. 
Moreover, as written in Remark 1.2, the appropriate truncations of unstable vortices remain unstable.
Let $\phi \in C^{\infty}_0 (B_1(0))$ be radially symmetric with $\phi \equiv 1$ on $B_{\frac{1}{2}}(0)$. 
For $\tilde R >0$, we set $\phi_{\tilde R} = \phi (\frac{\xi}{\tilde R})$ and the truncation of $\bar{u}$ and $\bar{\omega}$,
\[
    \bar{u}_{\tilde R} = \phi_{\tilde R} \bar{u}, 
    \quad
    \bar{\omega}_{\tilde R} = \nabla^{\bot} \cdot \bar{u}_{\tilde R},
\]
and also consider operators $\tilde \Lambda $ defined analogously to $\tilde \Lambda_{0}$ with 
$\bar{u}_{\tilde R}$ and $\bar{\omega}_{\tilde R}$
replacing $\bar{u}$ and $\bar{\omega}$.
\[
\begin{aligned}
    \tilde \Lambda  : 
    D(\tilde \Lambda ) \subset L^2_m (\mathbb{R}^2) \to L^2_m (\mathbb{R}^2), 
    &\quad 
    D(\tilde \Lambda ) 
    := 
    \{
        \omega \in L^2_m (\mathbb{R}^2) |\nabla \cdot (\bar{u}_{\tilde R} \omega) \in L^2_m (B_{\tilde R})
    \},
    \\
    \tilde \Lambda  \omega 
    &= 
    - \nabla \cdot (\bar{u}_{\tilde R} \omega)
    - (\tilde K [\omega] \cdot \nabla) \bar{\omega}_{\tilde R}.
\end{aligned}
\]
We now state the instability theorem by \cite{ABC1}.
\begin{prop} {\rm (\cite{ABC1})}
Let $\tilde \lambda$ be an unstable eigenvalue of $\tilde \Lambda_{0}$. 
For any $\epsilon \in (0, \Re \tilde \lambda)$, there exists $\tilde R_0 \ge 1$, such that for any $\tilde R \ge \tilde R_0$, the linear operator $\tilde \Lambda  : D(\tilde \Lambda ) \mapsto L^2_m (B_{\tilde R})$ has an unstable eigenvalue $\tilde \lambda_{\tilde R}$ with $|\tilde \lambda_{\tilde R} - \tilde \lambda| < \epsilon $.
\end{prop}
We note that $\tilde \lambda$ and $\tilde \lambda_{\tilde R}$ are isolated eigenvalues of $\tilde \Lambda_0$ and $\tilde \Lambda$, respectively.
We fix $R_0 \ge \tilde R_0$ and,
for simplicity, we set the unstable flow, the unstable vortex, and the isolated unstable eigenvalue as
\begin{equation}\label{061601}
    \tilde U^{st} := \bar{u}_{R_0}, 
    \quad 
    \tilde \Omega^{st} := 
    \bar{\omega}_{R_0},
    \quad 
    \lambda_{\infty} := 
    \tilde \lambda_{R_0}.
\end{equation}
We next check that $\lambda_{\infty}$ remains the unstable eigenvalue of the linear operator $\tilde \Lambda $ in $L^2 (B_{R_0}(0))$. 
For $\tilde U^{st}$ and $\tilde \Omega^{st}$,
the linear operator $\tilde \Lambda $ in $L^2 (B_{R_0}(0))$ is defined as 
\begin{equation}
   \begin{aligned} \label{061102}
      D (\tilde \Lambda )
      &= 
      \{
         \omega \in L^2 (B_{R_0}(0))
         |
         \nabla \cdot (\tilde U^{st} \omega ) \in L^2 (B_{R_0}(0))
      \},
      \\
      \tilde \Lambda  \omega
      &=
      - \nabla \cdot (\tilde U^{st}  \omega)
      - (\tilde K [ \omega \chi_{B_{R_0}(0)}] \cdot \nabla) \tilde \Omega^{st},
      \quad
     \omega \in D(\tilde \Lambda ).
   \end{aligned}
\end{equation}
As stated in the following proposition, it is easy to check that $\lambda_{\infty}$ remains the unstable eigenvalue of the linear operator $\tilde \Lambda $ in $L^2(B_{R_0}(0))$.
\begin{prop}\label{061201}
  The complex number $\lambda_{\infty}$ is an isolated unstable eigenvalue of the operator $\tilde \Lambda : D(\tilde \Lambda ) \subset L^2 (B_{R_0}(0)) \mapsto L^2 (B_{R_0}(0))$.
\end{prop}

\begin{proof}
   Let $\eta$ be a non-trivial unstable eigenfunction associated with an unstable eigenvalue $\lambda_{\infty}$. 
Since ${\rm supp\,} (\tilde \Lambda \eta) \subset B_{R_0}(0)$ and $\lambda_{\infty} \eta - \tilde \Lambda \eta = 0$, it follows that ${\rm supp\,} \eta \subset B_{R_0}(0)$.
In particular, $\eta \in D(\tilde \Lambda )$ and $\lambda_{\infty}$ is still an unstable eigenvalue of $\tilde \Lambda  : D(\tilde \Lambda ) \subset L^2 (B_{R_0}(0)) \mapsto L^2 (B_{R_0}(0))$.
\end{proof}

Using the profiles $\tilde U^{st}$, $\tilde \Omega^{st}$ and the linear operator $\tilde \Lambda $, we construct a velocity profile $\tilde U^{st}_{R}$ and a vorticity profile $\tilde \Omega^{st}_{R}$ with $R >R_0$, and a linear operator $\tilde \Lambda_{E,R}$ associated with $\tilde U^{st}_{R}$ and $\tilde \Omega^{st}_{R}$ in $L^2_{odd} (\mathbb{R}^2)$.
Here, $L^2_{odd} (\mathbb{R}^2)$ is given by
\[
\begin{aligned}
    L^2_{odd} (\mathbb{R}^2)
    := 
    \{
        \tilde g \in L^2 (\mathbb{R}^2)|
        \tilde g(\xi_1, -\xi_2)
        = 
        -\tilde g(\xi_1, \xi_2)
    \},
\end{aligned}
\]
with the standard $L^2$ norm. 
We first give the definition of the velocity profile $\tilde U^{st}_{R}$ and the vorticity profile $\tilde \Omega^{st}_{R}$ respectively by $\tilde U^{st}$ and $\tilde \Omega^{st}$. 
We define the velocity profile $\tilde U^{st}_{R}$ with $R >R_0$ as follows:
\begin{equation}
  \begin{aligned}\label{062901}
     \tilde U^{st}_{R} (\xi)
     = 
     \begin{pmatrix} 
       \tilde U^{st}_{R,1} (\xi) 
       \\  
       \tilde U^{st}_{R,2} (\xi)
    \end{pmatrix}
    &=
    \begin{pmatrix} 
       \tilde U^{st}_{1} (\xi_1, \xi_2 - R)
       \\
       \tilde U^{st}_{2} (\xi_1, \xi_2 - R)
    \end{pmatrix}+
       \begin{pmatrix}
        \tilde U^{st}_{1} (\xi_1, -\xi_2 - R) 
       \\  
        - \tilde U^{st}_{2} (\xi_1, -\xi_2 - R) 
    \end{pmatrix}
    \\
     & =: \tilde U^{st,+}_{R} (\xi) + \tilde U^{st,-}_{R} (\xi). 
 \end{aligned}
\end{equation}
We also give the vorticity profile $\tilde \Omega^{st}_{R} = \partial_1 \tilde U^{st}_{R,2} - \partial_2 \tilde U^{st}_{R,1}$, and then $\tilde \Omega^{st}_{R}$ can be written by
\[
\begin{aligned}
    \tilde \Omega^{st}_{R} (\xi)
    &= 
    \tilde \Omega^{st} (\xi_1, \xi_2 - R) - \tilde \Omega^{st} (\xi_1, - \xi_2 - R) 
    \\
    &=: \tilde \Omega^{st,+}_{R} (\xi) + \tilde \Omega^{st,-}_{R} (\xi). 
\end{aligned}
\]
Here, $\tilde \Omega^{st, \pm}_{R} = \partial_1 \tilde U^{st, \pm}_{R,2} - \partial_2 \tilde U^{st, \pm}_{R,1}$.
Let us give some remarks on $ \tilde U^{st}_{R}$ and $\tilde \Omega^{st}_{R}$.
For $\xi_2$ variable, 
$\tilde U^{st}_{R,1}$ is an even function and 
$\tilde U^{st}_{R,2}$ and $\tilde \Omega^{st}_{R}$ are odd functions, i.e.
$\tilde U^{st}_{R}$ and
$\tilde \Omega^{st}_{R}$ satisfy
\[
\begin{aligned}
      \tilde U^{st}_{R,1} (\xi_1, \xi_2)
      &= \tilde U^{st}_{R,1} (\xi_1, -\xi_2), & 
      \xi \in \mathbb{R}^2,
      \\
      \tilde U^{st}_{R,2} (\xi_1, \xi_2)
      &= -\tilde U^{st}_{R,2} (\xi_1, -\xi_2), & 
      \xi \in \mathbb{R}^2,
      \\
      \tilde \Omega^{st}_{R} (\xi_1, \xi_2)
      &= -\tilde \Omega^{st}_{R} (\xi_1, -\xi_2), &
      \xi \in \mathbb{R}^2. 
\end{aligned}
\]
Since the profiles $\tilde U^{st}$ and $\tilde \Omega^{st}$ have compact support in $\mathbb{R}^2$, and $\tilde U^{st}_{R}$ and $\tilde \Omega^{st}_{R}$ are their parallel translations in the $\xi_2$-direction, the supports of $\tilde U^{st}_{R}$ and $\tilde \Omega^{st}_{R}$ are also compact and are contained in the parallel translations of the support of the profiles:
\begin{equation} \label{061101}
   {\rm supp \,}\tilde U^{st, \pm}_{R}, \quad {\rm supp \,}\tilde \Omega^{st,\pm}_{R} 
   \subset 
   B_{R_0} ((0, \pm R)).
\end{equation}
We next define the linear operator $\tilde \Lambda_{E,R}$ associated with $\tilde U^{st}_{R}$ and $\tilde \Omega^{st}_{R}$ as follows:
\begin{equation}
    \begin{aligned} \label{061603}
        D (\tilde \Lambda_{E,R}) 
        &= 
        \{\tilde \omega \in 
        L^2_{odd}(B_{R_0}((0,R)) \cup B_{R_0}((0,-R))) 
           |
        \nabla \cdot (\tilde U^{st}_{R} \tilde \omega) \in L^2_{odd} (\mathbb{R}^2)
        \},
        \\
        \tilde \Lambda_{E,R} \tilde \omega (\xi)
        &:=
        - \nabla \cdot (\tilde U^{st}_{R} \tilde \omega) (\xi)
        - (\tilde K_R [\tilde \omega] \cdot \nabla)\tilde \Omega^{st}_{R} (\xi),
        \, \tilde\omega \in D (\tilde \Lambda_{E,R} ), 
    \end{aligned}
\end{equation}
where $\tilde K_R [\tilde \omega]$ and $\tilde \omega^{\pm}$ are given by
\[
\begin{aligned}
    \tilde K_R [\tilde \omega] 
    &= 
    \tilde K [\tilde \omega^{+}]
    +
    \tilde K [\tilde \omega^{-}],
    \\
    \tilde \omega^{+}
    &=
    \tilde \omega \chi_{B_{R_0}((0,R))},
    \quad
    \tilde \omega^{-}
    =
    \tilde \omega \chi_{B_{R_0}((0,-R))}.
\end{aligned}
\]
Here, $\tilde \omega$ and $\tilde K_R [\tilde \omega]$ satisfy the following properties. 
By the definition, $\tilde \omega$ and $\tilde \omega^{\pm}$ satisfy
\begin{equation}
   \begin{aligned} \label{061301}
      &\tilde \omega = \tilde \omega^{+} + \tilde \omega^{-},
      \\
      &{\rm supp\,}\tilde \omega^{+} \subset \overline{B_{R_0} ((0,R))}, 
      \quad
      {\rm supp\,}\tilde \omega^{-} \subset \overline{B_{R_0} ((0,-R))}.
   \end{aligned}
\end{equation}
We also note that the velocity $\tilde K_R [\tilde \omega]$ is a non-local function, in contrast to the support of $\tilde \omega$, $\tilde U^{st}_{R}$, and $\tilde \Omega^{st}_{R}$ are compact in $\mathbb{R}^2$. 
More precisely, 
$\tilde K_R [\tilde \omega]$ satisfies
\begin{equation}
\begin{aligned}\label{061602}
    \tilde K_R [\tilde \omega]_1 (\xi_1, \xi_2)
    &=
    \tilde K_R [\tilde \omega]_1 (\xi_1, -\xi_2), 
    \quad 
   & \xi \in \mathbb{R}^2,
    \\
    \tilde K_R [\tilde \omega]_2 (\xi_1, \xi_2)
    &=
    -\tilde K_R [\tilde \omega]_2 (\xi_1, -\xi_2), 
    \quad 
   & \xi \in \mathbb{R}^2. 
\end{aligned}
\end{equation}
Moreover, for any $\tilde \omega \in L^2_{odd}(B_{R_0}((0,R)) \cup B_{R_0}((0,-R)))$, there exists $\omega \in L^2 (B_{R_0}(0))$ such that $\tilde \omega^{+}$ and $\tilde \omega^{-}$ satisfy 
\[
    \begin{aligned}
        \tilde \omega^{+} (\xi_1, \xi_2)
        &=
         \omega (\xi_1, \xi_2 -R), 
        \quad 
       & \xi \in B_{R_0} ((0,R)),
        \\
        \tilde \omega^{-} (\xi_1, \xi_2)
        &=
        -\omega (\xi_1, -\xi_2 - R), 
        \quad 
       & \xi \in B_{R_0} ((0, -R)).
     \end{aligned}
\]
Since $\tilde U^{st}_{R}$, $\tilde \Omega^{st}_{R}$, and $\tilde \omega$ are compactly supported in $\mathbb{R}^2$ and the velocity $\tilde K_R [\tilde \omega]$ is a non-local function, 
we decompose $\tilde \Lambda_{E,R} = \tilde \Lambda^{+}_{E,R} + \tilde \Lambda^{-}_{E,R} + \tilde \Lambda^{remainder}_{E,R}$ with $\tilde \omega = \tilde \omega^{+} + \tilde \omega^{-}$ to
\[
\begin{cases}
    \tilde \Lambda^{+}_{E,R} \tilde \omega (\xi)
    &:= 
    -\nabla \cdot (\tilde U^{st,+}_{R} \tilde \omega^{+}) (\xi)
    -
    (\tilde K [\tilde \omega^{+}] \cdot \nabla) \tilde \Omega^{st,+}_{R} (\xi),
    \\
    \tilde \Lambda^{-}_{E,R} \tilde \omega (\xi)
    &:= 
    -
    \nabla \cdot (\tilde U^{st,-}_{R} \tilde \omega^{-}) (\xi)
    -
    (\tilde K [\tilde \omega^{-}] \cdot \nabla) \tilde \Omega^{st,-}_{R} (\xi),
    \\
    \tilde \Lambda^{remainder}_{E,R} \tilde \omega (\xi)
            &:= 
            -(\tilde K [\tilde \omega^{+}] \cdot \nabla) \tilde \Omega^{st,-}_{R} (\xi)
            -
            (\tilde K [\tilde \omega^{-}] \cdot \nabla) \tilde \Omega^{st,+}_{R} (\xi).
\end{cases}
\]
We note that $\nabla \cdot (\tilde U^{st,-}_{R} \tilde \omega^{+}) (\xi) = \nabla \cdot (\tilde U^{st,+}_{R}\tilde \omega^{-}) (\xi) = 0$, since the intersection of the supports $\tilde \omega^{+}$ and $\tilde U^{st,-}_{R}$, and the intersection of the supports $\tilde \omega^{-}$ and $\tilde U^{st,+}_{R}$ are empty. 
 By \eqref{061102} and \eqref{061603}, 
 $\tilde \Lambda^{\pm}_{E,R} : L^2 (B_{R_0} ((0,\pm R))) \mapsto L^2 (B_{R_0} ((0,\pm R)))$, and the following holds:
\begin{equation}
    \begin{aligned} \label{062701}
        \tilde \Lambda^{+}_{E,R} \tilde \omega (\xi_1, \xi_2)
        &=
        \tilde \Lambda   \omega (\xi_1, \xi_2 -R),
        \qquad 
        \xi \in B_{R_0}((0,R)),
        \\
        \tilde \Lambda^{-}_{E,R} \tilde \omega (\xi_1, \xi_2)
        &=
        -\tilde \Lambda  \omega (\xi_1, -\xi_2 -R),
        \quad
        \xi \in B_{R_0}((0,-R)).
     \end{aligned} 
\end{equation}
Finally, we demonstrate that the linear operator $\tilde \Lambda_{E,R}$ has an unstable eigenvalue.
To show the existence of an unstable eigenvalue for \eqref{061603}, we consider the resolvent problem of $\tilde \Lambda_{E,R}$. 
Let us fix $\lambda \in \rho_{re}(\tilde \Lambda)  \cap \{\Re \lambda >0\}$.
For given $\tilde g \in L^2_{odd}(B_{R_0}((0,R)) \cup B_{R_0}((0,-R)))$, we construct the solution $\tilde \omega \in L^2_{odd}(B_{R_0}((0,R)) \cup B_{R_0}((0,-R)))$ to the resolvent problem
\begin{equation}   \label{060201}
      (\lambda I - \tilde \Lambda_{E,R}) \tilde \omega = \tilde g,
\end{equation}
by taking a suitable approximation described as follows. 
Let $\tilde \omega^{\pm} [\tilde g^{\pm}]$ be defined as 
\begin{equation}
  \begin{cases} \label{062702}
      (\lambda - \tilde \Lambda^{+}_{E,R})\tilde \omega^{+} [\tilde g^{+}] (\xi)
      = \tilde g^{+} (\xi),
      \quad &\xi \in B_{R_0} ((0, R)), 
      \\
      (\lambda - \tilde \Lambda^{-}_{E,R})\tilde \omega^{-} [\tilde g^{-}] (\xi)
      = \tilde g^{-} (\xi),
      \quad &\xi \in B_{R_0} ((0, - R)), 
  \end{cases}
\end{equation}
where we denote $\tilde g =  \tilde g \chi_{B_{R_0}((0,R))} + \tilde g \chi_{B_{R_0}((0,-R))} =: \tilde g^{+} + \tilde g^{-} $. 
We remark that the supports of $\tilde \omega^{\pm} [\tilde g^{\pm}]$ are compact and satisfy 
${\rm supp\,}\tilde \omega^{\pm} [\tilde g^{\pm}] \subset \overline{B_{R_0}((0, \pm R))}$.
By using $\tilde \omega^{\pm}[\tilde g^{\pm}]$, we define an approximation $\tilde \omega_{app} [\tilde g]$ as follows:
\[
\begin{aligned}
    \tilde \omega_{app} [\tilde g]
    := 
    \tilde \omega^{+} [\tilde g^{+}] 
    + \tilde \omega^{-} [\tilde g^{-}].
\end{aligned}
\]
We comment that $\tilde \omega^{\pm}[\tilde g^{\pm}]$ can be bounded by $\tilde g$ with a constant which is independent of $R>R_0$. 
Indeed, when $g \in L^2(B_{R_0}(0))$ is such that 
\[
\begin{aligned}
    \tilde g^{+} (\xi_1, \xi_2)
    &=
    g (\xi_1, \xi_2 -R), 
    \quad 
    &\xi \in B_{R_0}((0,R)),
    \\
    \tilde g^{-}  (\xi_1, \xi_2)
    &=
    -g (\xi_1, -\xi_2 -R),
    \quad 
    &\xi \in B_{R_0}((0,-R)),
\end{aligned}
\]
we have by setting $\omega [g] = (\lambda I - \tilde \Lambda)^{-1} g$ in $L^2 (B_{R_0}(0))$, 
\[
\begin{aligned}
    \tilde \omega^{+} [\tilde g^{+}] (\xi_1, \xi_2)
    &= 
    \omega [g] (\xi_1, \xi_2 -R), 
    \quad 
    &\xi \in B_{R_0}((0,R)),
    \\
    \tilde \omega^{-} [\tilde g^{-}] (\xi_1, \xi_2) 
    &= 
    \omega [g] (\xi_1, -\xi_2 -R), 
    \quad 
    &\xi \in B_{R_0}((0,-R)).
\end{aligned}
\]
Since $\lambda \in \rho_{re} (\tilde \Lambda)$, we obtain
\begin{equation}
\begin{aligned} \label{062703}
    \| \tilde \omega^{\pm} [\tilde g]\|_{L^2(B_{R_0}(0, \pm R))}
    &= 
    \| (\lambda - \tilde \Lambda)^{-1} g\|_{L^2(B_{R_0}(0))}
    \\
    &\le
    C \| g\|_{L^2(B_{R_0}(0))}
    =
    C \| \tilde g^{\pm} \|_{L^2  (B_{R_0}((0, \pm R)))}.
\end{aligned}
\end{equation}
Here, the constant $C$ in the right-hand side depends on $R_0$ and the operator norm of $(\lambda -\tilde \Lambda)^{-1}$, and 
the $L^2$ norm of $\tilde g$ does not depend on $R > R_0$.
Therefore, \eqref{062703} implies that $\tilde \omega^{\pm}[\tilde g^{\pm}]$ is estimated independently $R> R_0$.
Using $\tilde \omega_{app}[\tilde g]$, we can give a formula which represents the solution $\tilde \omega$. 
Acting the remainder operator $\tilde \Lambda^{remainder}_{E,R}$ on $\tilde \omega_{app} [\tilde g]$, we obtain
\begin{equation}
    \begin{aligned} \label{052902}
        \tilde \Lambda^{remainder}_{E,R} \tilde \omega_{app}[\tilde g] (\xi)
        &= 
        (\tilde K [\tilde \omega^{+} [\tilde g^{+}] ] \cdot \nabla) \tilde \Omega^{st,-}_{R} (\xi)
        \\
        &\hspace{11pt}+
        (\tilde K [\tilde \omega^{-} [\tilde g^{-}] ] \cdot \nabla) \tilde \Omega^{st,+}_{R} (\xi)
        \\
        &=: U_R \tilde g (\xi).
    \end{aligned}
\end{equation}
Combining the definition of $\tilde \omega_{app}[\tilde g]$ and \eqref{052902} leads to the following modified equation:
\[
\begin{aligned}
    (\lambda- \tilde \Lambda_{E,R}) 
    \tilde \omega_{app} [\tilde g]
    =
    (I - U_R) \tilde g.
\end{aligned}
\]
In the next lemma, we give the estimate of the remainder term $U_R \tilde g$.
\begin{lem}{\rm (Estimate of $U_R$)}
    \label{052302}
For any $\tilde g \in L^2_{odd}(B_{R_0}((0,R)) \cup B_{R_0}((0,-R))) $, 
\begin{equation} \label{052901}
    \| U_R \tilde g \|_{L^2 (\mathbb{R}^2)} 
    \le
    \frac{C}{R - R_0} \|\tilde g\|_{L^2 (\mathbb{R}^2)},
\end{equation}
where $C$ depends on $R_0$, $\tilde \Omega^{st}$, and the operator norm of $(\lambda - \tilde \Lambda)^{-1}$ on $L^2 (B_{R_0}(0))$.
\end{lem}

\begin{proof}
    By symmetry, it suffices to show the estimate for $\tilde K [\tilde \omega^{+}[\tilde g^{+}]]$ in \eqref{052902}. 
    By the Biot--Savart law, $ \tilde K [\tilde \omega^{+}[\tilde g^{+}]] = \nabla^{\bot} (- \Delta)^{-1} \tilde \omega^{+}[\tilde g^{+}]$ can be written as
\[
\begin{aligned}
    \tilde K [\tilde \omega^{+}[\tilde g^{+}]] (\xi) := 
    \frac{1}{2 \pi}
    \int_{B_{R_0}((0,R))} 
    \frac{(\xi -\eta)^{\bot}}{|\xi -\eta|^2}
    \tilde \omega^{+}[\tilde g^{+}] (\eta) d\eta.
\end{aligned}
\]
   By \eqref{061101} and \eqref{061301}, $\tilde \omega^{+}[\tilde g^{+}] (\eta)$ has a compact support in $\mathbb{R}^2_+$ and is located far from the support of $\tilde \Omega^{st, -}_{R}(\xi)$, in particular, $|\xi -\eta| \ge 2(R - R_0)$. 
  Using the estiamte for $\tilde K [\tilde \omega^{+}[\tilde g^{+}]]$, Fubini's theorem, and H\"older's inequality, we obtain
\[
\begin{aligned}
    &\| (\tilde K [\tilde \omega^{+}[\tilde g^{+}]] \cdot \nabla) \tilde \Omega^{st,-}_{R}
    \|_{L^2(B_{R_0}((0, -R)))}
    \\
    \le & 
    \frac{C}{R-R_0}     
    \|\nabla \tilde \Omega^{st,-}_{R}(\xi)\|_{L^2(B_{R_0}((0, -R)))}
    \int_{B_{R_0}((0,R))}
    |\tilde \omega^{+}[\tilde g^{+}]|
            d\eta 
            \\
            \le&
            \frac{C(R_0, \tilde \Omega^{st})}{R-R_0} 
            \|\tilde \omega^{+}[\tilde g^{+}]\|_{L^2(B_{R_0}((0, R)))}.
\end{aligned}
\]
 Since $\lambda \in \rho_{re} (\tilde \Lambda)$ and \eqref{062703},
 hence
\[
\begin{aligned}
    \| (\tilde K[\tilde \omega^{+} [\tilde g^{+}]] \cdot \nabla) \tilde \Omega^{st,-}_{R}
    \|_{L^2(B_{R_0}((0,-R)))}
    \le 
    \frac{C}{R-R_0} \|\tilde g^{+}\|_{L^2 (B_{R_0}((0,R)))}.
\end{aligned}
\]
The proof is complete. 
\end{proof}
Since the constant in Lemma \ref{052302} relies on $R_0$, $\tilde \Omega^{st}$, and the operator norm of the resolvent $(\lambda -\tilde \Lambda)^{-1}$,
  the operator $I - U_R$ is invertible in $L^2_{odd} (B_{R_0} ((0,R)) \cup B_{R_0} ((0,-R)))$ for sufficiently large $R > C(R_0, \tilde \Omega^{st}, (\lambda -\tilde \Lambda)^{-1})$. 
In particular, we take $(I - U_R)^{-1}\tilde g$ in stead of $\tilde g$, then 
\[
    (\lambda - \tilde \Lambda_{E,R}) 
    \tilde \omega_{app} [(I - U_R)^{-1}\tilde g]
    = 
    \tilde g.
\]
By the above equation, we obtain the representation of the solution $\tilde \omega$ to \eqref{060201} by $\tilde \omega_{app}$ as follows:
\[
  \tilde \omega
    = \tilde \omega_{app}[(I -U_R)^{-1}\tilde g].
\]
Moreover, using the Neumann series, we have
\begin{equation}
    \begin{aligned} \label{062302}
        \tilde \omega 
        = 
        \tilde \omega_{app} [\tilde g]
        + 
        \tilde \omega_{app} \bigg[\sum_{k \in \mathbb{N}}  U_R^k \tilde g \bigg].
    \end{aligned} 
\end{equation}

\begin{prop} \label{052301}
    Let $\lambda_{\infty}$ be an isolated unstable eigenvalue of $\tilde \Lambda $. 
    For any $\epsilon \in (0, {\Re\,}\lambda_{\infty})$, there exists $R' > R_0$ which depends on $\epsilon, \lambda_{\infty}, R_0$, and $\tilde \Omega^{st}$
     such that for any $R \ge R'$, the operator $\tilde \Lambda_{E,R}: D(\tilde \Lambda_{E,R}) \mapsto L^2_{odd}(B_{R_0}((0,R))\cup B_{R_0}((0,-R)))$ has an unstable eigenvalue $\lambda_{R}$ with $\lambda_{R} \in B_{\epsilon}(\lambda_{\infty})$. 
    \end{prop}

\begin{proof}
  Let $\vec{c} = \partial B_{\epsilon} (\lambda_{\infty})$. 
  Since $\lambda_{\infty}$ is an isolated unstable eigenvalue of $\tilde \Lambda $, we take $\epsilon > 0$ so small that $\sigma (\tilde \Lambda) \cap (\overline{B_{\epsilon}(\lambda_{\infty})} \setminus \{\lambda_{\infty} \}) = \emptyset$, and hence we define the spectral projection $Pr : L^2 (B_{R_0}(0))\mapsto L^2 (B_{R_0}(0))$ with $\tilde \Lambda$ as
  \[
  \begin{aligned}
  Pr 
  = 
  \frac{1}{2 \pi i} \int_{\vec{c}} 
  (\lambda - \tilde \Lambda)^{-1} d\lambda.
  \end{aligned}
  \]
  In particular, $Pr$ is non-trivial. 
  We observe from 
  \[
  \sup_{\lambda \in \vec{c}} \|(\lambda - \tilde \Lambda)^{-1}\|_{L^2 (B_{R_0}(0)) \to L^2 (B_{R_0}(0))}< \infty,
  \]
 and from the construction of the solution \eqref{062302} to the resolvent problem \eqref{060201} 
 that $\vec{c} \subset \rho_{re}(\tilde \Lambda_{E,R})$ for sufficiently large $R$.
  Hence we can also define the spectral projection $Pr_{R} : L^2_{odd}(B_{R_0}((0,R)) \cup B_{R_0}((0,-R)))  \mapsto L^2_{odd}(B_{R_0}((0,R)) \cup B_{R_0}((0,-R)))$ such that
  \[
  \begin{aligned}
  Pr_{R} 
  = 
  \frac{1}{2\pi i} \int_{\vec{c}} 
  (\lambda I - \tilde \Lambda_{E,R})^{-1} d\lambda.
  \end{aligned}
  \]
  Our aim is to show that $Pr_R$ is also non-trivial. 
  Let $\tilde h \in L^2 (B_{R_0}(0))$ be an eigenfunction for the isolated unsable eigenvalue $\lambda_{\infty}$ of $\tilde \Lambda$. 
  We set $\tilde g = \tilde g^{+} + \tilde g^{-}$ by $\tilde h$ as follows:
  \[
    \begin{cases}
  \tilde g^{+} (\xi_1, \xi_2)
  = \tilde h\chi_{B_{R_0}(0)} (\xi_1, \xi_2 -R), 
  \qquad 
         &\xi \in B_{R_0} ((0,R)),
         \\
  \tilde g^{-} (\xi_1, \xi_2)
  = -\tilde h\chi_{B_{R_0}(0)} (\xi_1, -\xi_2 - R), 
  \quad 
         &\xi \in B_{R_0} ((0,-R)). 
    \end{cases}
    \]
  Then we have $\tilde g \in L^2_{odd}(B_{R_0}((0,R)) \cup B_{R_0}((0,-R)))$, and for $\lambda \in \rho_{re} (\tilde \Lambda)$,  $ \tilde \omega_{app}[\tilde g]$ satisfies
  \[
  \tilde \omega_{app} [\tilde g]
  = 
  \tilde \omega^{+} [\tilde g^{+}]
  +
  \tilde \omega^{-} [\tilde g^{-}]
  = 
  (\lambda I - \tilde \Lambda^{+}_{E,R})^{-1} \tilde g^{+}
  +
  (\lambda I - \tilde \Lambda^{-}_{E,R})^{-1} \tilde g^{-}.
  \]
  In particular, combining the fact that $\tilde \Lambda^{\pm}_{E,R}$ are the parallel $\pm R$ shift of $\tilde \Lambda$ and the above equation, the following holds:
  \begin{equation}
      \begin{cases} \label{062801}
  (\lambda I - \tilde \Lambda^{+}_{E,R})^{-1} \tilde g^{+} (\xi_1, \xi_2)
  =
  (\lambda I - \tilde \Lambda)^{-1} \tilde h (\xi_1, \xi_2 -R), 
          \\
  (\lambda I - \tilde \Lambda^{-}_{E,R})^{-1} \tilde g^{-} (\xi_1, \xi_2)
  =
  -(\lambda I - \tilde \Lambda)^{-1} \tilde h (\xi_1, -\xi_2 - R). 
      \end{cases}
  \end{equation}
  By using \eqref{062302} and the spectral projection, we obtain
  \begin{equation}
      \begin{aligned} \label{062802}
          &\frac{1}{2\pi i} \int_{\vec{c}} 
  (\lambda I - \tilde \Lambda_{E,R})^{-1} 
  \tilde g 
  d\lambda 
          \\
  =
          &
  \frac{1}{2\pi i} \int_{\vec{c}} 
  (\lambda I - \tilde \Lambda^{+}_{E,R})^{-1} \tilde g^{+} 
  d\lambda 
  +
  \frac{1}{2\pi i} \int_{\vec{c}} 
  (\lambda I - \tilde \Lambda^{-}_{E,R})^{-1} \tilde g^{-}
  d\lambda 
          \\
          &
  + \frac{1}{2\pi i} \int_{\vec{c}} 
  (\lambda I - \tilde \Lambda^{+}_{E,R})^{-1} 
           \sum_{k \in \mathbb{N}} (U_R^k \tilde g)^{+} 
  d\lambda 
  +
  \frac{1}{2\pi i} \int_{\vec{c}} 
  (\lambda I - \tilde \Lambda^{-}_{E,R})^{-1} 
          \sum_{k \in \mathbb{N}} (U_R^k \tilde g)^{-}
  d\lambda 
          \\
  =&
  \bigg(
  \frac{1}{2\pi i} \int_{\vec{c}} 
  (\lambda I - \tilde \Lambda )^{-1} 
  \tilde h 
  d\lambda 
  \bigg) (\xi_1, \xi_2 -R)
  -
  \bigg(
  \frac{1}{2\pi i} \int_{\vec{c}} 
  (\lambda I - \tilde \Lambda )^{-1} 
  \tilde h 
  d\lambda 
  \bigg) (\xi_1, -\xi_2 - R)
          \\
          & +
  \frac{1}{2\pi i} \int_{\vec{c}} 
  (\lambda I - \tilde \Lambda^{+}_{E,R})^{-1}
           \sum_{k \in \mathbb{N}} (U_R^k \tilde g)^{+} 
  d\lambda 
  +
  \frac{1}{2\pi i} \int_{\vec{c}} 
  (\lambda I - \tilde \Lambda^{-}_{E,R})^{-1}
           \sum_{k \in \mathbb{N}} (U_R^k \tilde g)^{-}
  d\lambda
           \\
  =&
  \tilde h (\xi_1, \xi_2 -R)
  -
  \tilde h (\xi_1, -\xi_2 -R)
           \\
           & +
  \frac{1}{2\pi i} \int_{\vec{c}} 
  (\lambda I - \tilde \Lambda^{+}_{E,R})^{-1}
            \sum_{k \in \mathbb{N}} (U_R^k \tilde g)^{+} 
  d\lambda 
  +
  \frac{1}{2\pi i} \int_{\vec{c}} 
  (\lambda I - \tilde \Lambda^{-}_{E,R})^{-1}
            \sum_{k \in \mathbb{N}} (U_R^k \tilde g)^{-}
  d\lambda
  .
      \end{aligned}
  \end{equation}
  Here, we have used $Pr \tilde h = \tilde h$ in the first and second terms of the right-hand side above, since $Pr$ is the spectral projection of $\tilde \Lambda$ and $\tilde h$ is an eigenfunction of $\tilde \Lambda$.
  We first focus on the first and second terms of the right-hand side in \eqref{062802}.
  Considering the support of $\tilde h$, the supports of the first and second terms of the right-hand side satisfy
  \[
  \begin{aligned}
  {\rm supp \,}  \tilde h  (\xi_1, \xi_2 -R) 
      \subset B_{R_0} ((0,R)), 
      \\
  {\rm supp \,} \tilde h (\xi_1, - \xi_2 - R) 
      \subset B_{R_0} ((0, - R)).
  \end{aligned}
  \]
  In particular, we remark that $\tilde h  (\xi_1, \xi_2 -R)$ does not cancel by $-\tilde h (\xi_1, - \xi_2 - R)$.
 Considering the $L^2$ norm of the first and second terms, we obtain
 \[
 \begin{aligned}
     &\| \tilde h (\xi_1, \xi_2 -R) 
  - \tilde h (\xi_1, -\xi_2 -R) \|^2_{L^2(B_{R_0}((0,R)) \cup B_{R_0}((0,-R)))}
     \\ 
  =& 
     \| \tilde h (\xi_1, \xi_2 -R) \|^2_{L^2(B_{R_0}((0,R)))}
  + 
     \| \tilde h (\xi_1, -\xi_2 -R) \|^2_{L^2(B_{R_0}((0,-R)))}
  =    2 \| \tilde h \|^2_{L^2(B_{R_0}(0))},
 \end{aligned}
 \]
 for $R>R_0$. 
 In particular, the $L^2$ norm of the first and second terms is independent of $R >R_0$. 
 We next focus on the third and fourth terms in \eqref{062802}.
  By Lemma \ref{052302}, the third and fourth terms in \eqref{062802} converge to $0$ in $L^2 (B_{R_0}((0,R)) \cup B_{R_0}((0,-R)))$ as $R \to \infty$ uniformly in $\lambda \in \vec{c}$.
  Then, for any $ \epsilon' >0$, there exists $R' > 0$, for any $R \ge R'$, the third term of the right-hand side in \eqref{062802} satisfies
  \[
  \begin{aligned}
  \bigg\|
  \frac{1}{2\pi i} \int_{\vec{c}} 
  (\lambda I - \tilde \Lambda^{+}_{E,R})^{-1} 
      \sum_{k \in \mathbb{N}} (U_R^k \tilde g)^{+}d\lambda 
  \bigg\|_{L^2 (B_{R_0}((0,R)))}
  < \frac{\epsilon'}{2},
  \end{aligned}
  \]
  and the fourth term also does.
 Therefore, $Pr_R$ is non-trivial for sufficiently large $R \ge R'$. 
  The proof is complete.
  \end{proof}
  We have considered the linear instability of the stationary velocity and vorticity profiles that are symmetric about the $\xi_1$-axis. 
  By \eqref{061602}, in particular, the vorticity $\tilde \omega$ and the normal component of the velocity $\tilde K_R [\tilde \omega]_2$ are $0$ at $\xi_2 = 0$. 
  For sufficiently large $R$ satisfying the assumption in Proposition~\ref{052301}, define $U_E(\xi) := \tilde U^{st, +}_R(\xi)$, which is defined in \eqref{062901}, and let $\Lambda_E$ be the operator in $L^2(B_{R_0}((0,R)))$ given by
  \[
  \begin{aligned}
    D(\Lambda_E) & =\{\omega \in L^2 (B_{R_0}((0,R))) | \nabla \cdot ( U_E \omega) \in L^2 (B_{R_0}((0,R)))\}\\
    \Lambda_E \omega & = - \nabla \cdot (U_E \omega) - K[\omega \chi_{B_{R_0}(0,R)}] \cdot \nabla \Omega_E,
    \quad 
    \omega \in D(\Lambda_E),
  \end{aligned}
  \]
  where $\Omega_E (\xi) = \tilde \Omega^{st,+}_R (\xi)$, and $K[\omega \chi_{B_{R_0}(0,R)}]$ denotes the velocity in $\mathbb{R}^2_{+}$ defined by \eqref{def.K}.
  Then $\Lambda_E$ in $L^2 (B_{R_0}((0,R)))$ and 
   $\tilde \Lambda_{E,R}$ in $L^2 (B_{R_0}((0,R))\cup B_{R_0}((0,-R)))$ are equivalent in the sense that the following relation holds:
   \[
   \begin{aligned}
       \Lambda_E
       = \gamma_{\mathbb{R}^2_{+}} \tilde \Lambda_{E,R} e_{odd},
       \\
       \tilde \Lambda_{E,R}
       = e_{odd} \Lambda_E \gamma_{\mathbb{R}^2_{+}}.
    \end{aligned}
   \]
   Here, $\gamma_{\mathbb{R}^2_{+}}$ denotes the restriction on $\mathbb{R}^2_{+}$ and $e_{odd}$ denotes the odd extension. 
   Therefore, the existence of an isolated unstable eigenvalue of $\tilde \Lambda_{E,R}$ implies the existence of an isolated unstable eigenvalue of $\Lambda_E$, in particular, we have Theorem \ref{thm.euler}.

\appendix

\section{Proof of Lemma \ref{lem.bs}}

To show \eqref{est.lem.bs.1} we consider the Poisson equations for the streamfunction $\varphi$, that is, $-\Delta \varphi = \omega$ with $\varphi|_{\xi_2=0} =0$. Then we have from the integration by parts,
\begin{align*}
\| \nabla \varphi \|^2_{L^2(\R^2_+)} = \langle -\Delta \varphi, \varphi \rangle_{L^2(\R^2_+)} &= \langle \omega,\varphi \rangle_{L^2(\R^2_+)}\\
& \leq \|\xi_2 \omega \|_{L^2(\R^2_+)} \| \frac{\varphi}{\xi_2}\|_{L^2(\R^2_+)}\\
& \leq C \|\xi_2 \omega \|_{L^2(\R^2_+)} \| \partial_2 \varphi \|_{L^2(\R^2_+)}\,,
\end{align*}
where the Hardy inequality is used in the last line. 
The estimate $\|\nabla^2 \varphi\|_{L^2(\R^2_+)}\leq \| \omega \|_{L^2(\R^2_+)}$ is a standard elliptic estimate in $L^2$ and is straightforward from the integration by parts. The last estimate in \eqref{est.lem.bs.1} follows from the interpoation inequality $\|\partial_2\varphi\|_{L^\infty(\R^2_+)}^2\leq C \| \partial_2^2 \varphi \|_{L^2(\R^2_+)} \| \partial_2 \varphi \|_{L^2(\R^2_+)}$ and from the other two estimates in \eqref{est.lem.bs.1}. As for \eqref{est.lem.bs.2} and \eqref{est.lem.bs.3}, we give a proof only for \eqref{est.lem.bs.3}, since \eqref{est.lem.bs.2} is proved  in a similar way. 
Let $\chi_\kappa(\xi_2)$, $\kappa>0$, be a cut-off such that $\chi_\kappa=1$ for $0\leq \xi_2\leq \kappa$ and $\chi_\kappa=0$ for $\xi_2\geq 2\kappa$. We note that $\xi_1\partial_1 \varphi$ satisfies 
\begin{align*}
-\Delta (\xi_1\partial_1\varphi) = -2 \partial_1^2 \varphi + \xi_1\partial_1\omega\,,
\end{align*}
which gives 
\begin{align*}
\| \nabla (\xi_1\partial_1\varphi)\|_{L^2(\R^2_+)}^2 & = \langle -2 \partial_1^2 \varphi + \xi_1\partial_1\omega, \xi_1\partial_1\varphi\rangle_{L^2(\R^2_+)} \\
& = \langle 2\partial_1\varphi - \xi_1\omega, \partial_1 (\xi_1\partial_1\varphi)\rangle_{L^2(\R^2_+)}-\langle \omega, \xi_1\partial_1\varphi\rangle_{L^2(\R^2_+)}\,.
\end{align*}
Thus, we obtain 
\begin{align}\label{proof.lem.bs.1}
\| \nabla (\xi_1\partial_1\varphi) \|_{L^2(\R^2_+)} &\leq C(\| \partial_1\varphi\|_{L^2(\R^2)} + \| \xi_1\omega\|_{L^2(\R^2_+)} + \| \xi_2\omega\|_{L^2(\R^2_+)}) \nonumber \\
&  \leq C\| |\xi| \omega\|_{L^2(\R^2_+)}\,.
\end{align}
Since $\chi_\kappa\xi_1\partial_1\varphi$ satisfies 
\begin{align}\label{proof.lem.bs.2}
-\Delta ( \chi_\kappa\xi_1\partial_1\varphi ) = -\chi_\kappa''\xi_1\partial_1\varphi -2\chi_\kappa' \xi_1 \partial_{12}^2 \varphi -2 \chi_\kappa \partial_1^2\varphi + \chi_\kappa \xi_1\partial_1\omega\,.
\end{align}
we have from $-\chi_\kappa'' \xi_1\partial_1\varphi = -\chi_\kappa''\int_0^{\xi_2} \partial_2 (\xi_1\partial_1\varphi ) \, d\eta_2$, 
\begin{align}\label{proof.lem.bs.3}
& \| \nabla^2 (\chi_\kappa \xi_1\partial_1\varphi) \|_{L^2(\R^2_+)} \nonumber \\
& \leq C\big ( \|\chi_\kappa''\xi_1\partial_1\varphi\|_{L^2(\R^2_+)}  + \| \chi_\kappa' \xi_1 \partial_{12}^2 \varphi\|_{L^2(\R^2_+)} + \|\chi_\kappa \partial_1^2\varphi\|_{L^2(\R^2_+)} + \| \chi_\kappa \xi_1\partial_1\omega\|_{L^2(\R^2_+)}\big ) \nonumber \\
& \leq C \big (\| \partial_2 (\xi_1\partial_1\varphi)\|_{L^2(\R^2_+)} + \| \partial_1^2 \varphi\|_{L^2(\R^2_+)} + \| \chi_{\kappa}\xi_1\partial_1\omega \|_{L^2(\R^2_+)}\big ) \nonumber \\
& \leq C \big ( \| \langle \xi\rangle \omega\|_{L^2(\R^2_+)} + \| \chi_{\kappa} \xi_1\partial_1\omega \|_{L^2(\R^2_+)}\big )\,.
\end{align}
Here we have used \eqref{est.lem.bs.1} and \eqref{proof.lem.bs.1} in the last line.
Hence we have from $\partial_2(\chi_\kappa \xi_1\partial_1\varphi) = \chi_\kappa \partial_2 (\xi_1\partial_1\varphi) + \chi_\kappa' \int_0^{\xi_2} \partial_2 (\xi\partial_1\varphi) d\eta_2$, 
\begin{align}
\| \xi_1\partial_{12}^2\varphi|_{\xi_2=0}\|_{L^2(\R)} & = \| \{\partial_2( \chi_\kappa \xi_1\partial_1\varphi) \}|_{\xi_2=0} \|_{L^2(\R)} \nonumber \\
& \leq C\| \partial_2^2  ( \chi_\kappa \xi_1\partial_1\varphi) \|_{L^2(\R^2_+)}^\frac12 \| \partial_2 ( \chi_\kappa \xi_1\partial_1\varphi) \|_{L^2(\R^2_+)}^\frac12\nonumber \\
& \leq C\| \partial_2^2  ( \chi_\kappa \xi_1\partial_1\varphi) \|_{L^2(\R^2_+)}^\frac12 \| \partial_2  (\xi_1\partial_1\varphi) \|_{L^2(\R^2_+)}^\frac12 \nonumber \\
& \leq C (\|\langle \xi\rangle\omega\|_{L^2(\R^2_+)} + \| \chi_\kappa \xi_1\partial_1\omega\|_{L^2(\R^2_+)})\,.
\end{align}
The case $j=1$ of \eqref{est.lem.bs.3} is obtained from the equations for $\chi_\kappa\xi_1\partial_1^2\varphi$  similarly by using \eqref{est.lem.bs.3} for $j=0$ and \eqref{est.lem.bs.1}, \eqref{est.lem.bs.2}. We omit the details. The proof is complepte.

\section{Proof of Lemma \ref{lem.stokes}}

We note that the integration by parts implies $\Re \langle \xi\cdot \nabla v+ v, v\rangle_{L^2(\R^2_+)}=0$. Hence, we have for $v=(\lambda I - \mathbb{H}_\alpha)^{-1}f$ the energy equality $\Re(\lambda) \| v\|_{L^2(\R^2_+)}^2 + \frac{1}{\alpha} \| \nabla v\|_{L^2 (\R^2_+)}^2 = \Re \langle f, v\rangle_{L^2(\R^2_+)}$ gives the esitmates for $\| v\|_{L^2(\R^2_+)}$ and $\| \nabla v\|_{L^2(\R^2_+)}$.
To show the estimate for $\nabla^2 v$, we observe that $\Re \langle (\lambda I - \mathbb{H}_\alpha) v, -\mathbb{A} v\rangle_{L^2(\R^2_+)} = \Re \langle f,-\mathbb{A}v\rangle_{L^2(\R^2_+)}$ gives 
\begin{align*}
& \Re (\lambda) \| \nabla  v\|_{L^2(\R^2_+)}^2 +\frac1\alpha \| \mathbb{A} v\|_{L^2(\R^2_+)}^2 - \frac1\alpha \Re \langle (\frac{\xi}{2}\cdot \nabla +\frac12) v, -\mathbb{A} v\rangle_{L^2(\R^2_+)} \\
& \quad = \Re \langle f,-\mathbb{A}v\rangle_{L^2(\R^2_+)}\,.
\end{align*}
Since $\frac{\xi}{2}\cdot \nabla v, \, v \in L^2_\sigma (\R^2_+)$, we see from the integration by parts,
\begin{align*}
\Re \langle (\frac{\xi}{2}\cdot \nabla +\frac12) v, -\mathbb{A} v\rangle_{L^2(\R^2_+)} & =  \Re \langle (\frac{\xi}{2}\cdot \nabla +\frac12) v, -\Delta v\rangle_{L^2(\R^2_+)} \\
& =  \Re \langle \frac{\xi}{2}\cdot \nabla \nabla v, \nabla v\rangle_{L^2(\R^2_+)} + \| \nabla v\|_{L^2(\R^2_+)}^2 \\
& = \frac12 \| \nabla v\|_{L^2(\R^2_+)}^2\,.
\end{align*}
Thus we have 
\begin{align*}
\frac1\alpha \| \mathbb{A} v\|_{L^2(\R^2_+)} \leq C \| f\|_{L^2(\R^2_+)}\,,
\end{align*}
which also gives $\frac1\alpha \| \nabla^2 v\|_{L^2(\R^2_+)}\leq  C \| f\|_{L^2 (\R^2_+)}$ by the elliptic estimate of the Stokes operator in the half space $\| \nabla^2 v\|_{L^2(\R^2_+)}\leq C \| \mathbb{A} v\|_{L^2(\R^2_+)}$, as desired. Then $\frac1\alpha \|\xi\cdot \nabla v\|_{L^2(\R^2_+)}\leq C\|f\|_{L^2(\R^2_+)}$ follows from the equations. The proof of \eqref{est.lem.stokes.1} is complete.
To show \eqref{est.lem.stokes.2} we see that the vorticity $\omega ={\rm rot}\, v$ satisfies 
\begin{align}
\lambda \nabla^k \omega -\frac{1}{\alpha} (\Delta +\frac{\xi}{2}\cdot \nabla + \frac{2+k}{2} ) \nabla^k \omega = \nabla^k{\rm rot}\, f\,,\quad k=0,1\,.
\end{align}
For a smooth bounded domain $O_k$ satisfying $\overline{O_k}\subset \R^2_+$, let $\chi_{O_k}$ be a cut-off function with compact support such that $\chi_{O_k}=1$ on $O_k$ and ${\rm supp}\, \chi_{O_k}\subset \R^2_+$. 
Then $\omega_k = \nabla^k \omega \, \chi_{O_k}$ satisfies 
\begin{align*}
& \lambda \omega_k -\frac{1}{\alpha} (\Delta +\frac{\xi}{2}\cdot \nabla + \frac{2+k}{2} ) \omega_k \\
& = -\frac1\alpha \big ((\Delta \chi_{O_k})\nabla^k \omega + 2\nabla \chi_{O_k}\cdot \nabla \nabla^k \omega + (\frac{\xi}{2}\cdot \nabla \chi_{O_k}) \nabla^k \omega \big )  + \chi_{O_k}\nabla^k  {\rm rot}\, f\,.
\end{align*}
Note that $\omega_k|_{\xi_2=0}=0$. Then, as in the estimate of $v$, we have 
\begin{align*}
& \|\omega_0\|_{L^2(\R^2_+)} + \frac{1}{\alpha^\frac12}\| \nabla \omega_0\|_{L^2(\R^2_+)} + \frac1\alpha \| \nabla^2 \omega_0\|_{L^2(\R^2_+)} \\
& \leq \frac{C}{\alpha} \sum_{j=0,1}\| \nabla^j\omega \|_{L^2 (\R^2_+)} + C \| \chi_{O_0} {\rm rot}\, f\|_{L^2(\R^2_+)}\\
& \leq C \|f\|_{L^2(\R^2_+)} +  C \| \chi_{O_0} {\rm rot}\, f\|_{L^2(\R^2_+)} \,,
\end{align*}
and then, by taking $O_0$ and $O_1$ with $\overline{O_1}\subset O_0$ suitablly, we have 
\begin{align*}
\|\omega_1\|_{L^2(\R^2_+)} + \frac1\alpha \| \nabla^2 \omega_1\|_{L^2(\R^2_+)} & \leq \frac{C}{\alpha} \sum_{j=0,1,2}\| \nabla^j\omega_0 \|_{L^2 (\R^2_+)} + C \| \chi_{O_1} \nabla {\rm rot}\, f\|_{L^2(\R^2_+)}\\
& \leq C \|f\|_{L^2(\R^2_+)} +  C\sum_{j=0,1} \| \chi_{O_0} \nabla^j {\rm rot}\, f\|_{L^2(\R^2_+)}\,.
\end{align*}
The proof is complete.

\

\noindent \textbf{Acknowledgements:}  The first author acknowledges the support of JSPS KAKENHI Grant Number JP24KJ0116, JP25K17276.
The second author acknowledges the support of JSPS KAKENHI Grant Number 21H04433, 23H01082, 23K25779, 24H00185, 24H00184.

\end{document}